%% file: DD_OpInf_final.tex
\documentclass[conf]{new-aiaa}

\input{settings}

\title{Domain decomposition for data-driven reduced modeling of large-scale systems}

\author{Ionu\cb{t}-Gabriel Farca\cb{s} \footnote{Postdoctoral fellow, Oden Institute for Computational Engineering and Sciences, ionut.farcas@austin.utexas.edu, AIAA Member.}}
\affil{The University of Texas at Austin, Austin, TX, 78712}
\author{Rayomand P. Gundevia  \footnote{Aerospace Research Engineer, rayomand.gundevia.ctr@afrl.af.mil, AIAA Member.}}
\affil{Jacobs Engineering Group, Inc., Edwards Air Force Base, CA 93524}
\author{Ramakanth Munipalli\footnote{Senior Aerospace Research Engineer, Combustion Devices, ramakanth.munipalli@us.af.mil, AIAA Senior Member.}}
\affil{Air Force Research Laboratory, Edwards Air Force Base, CA 93524}
\author{Karen E.~Willcox\footnote{Director, Oden Institute for Computational Engineering and Sciences, kwillcox@oden.utexas.edu, AIAA Fellow.}}
\affil{The University of Texas at Austin, Austin, TX, 78712}

\begin{document}

\maketitle
\begin{tikzpicture}[remember picture,overlay]
    \node[anchor=base west] at (current page footer area.south west) {Distribution Statement A: Approved for Public Release; Distribution is Unlimited. PA$\#$ AFRL-2023-5401};
\end{tikzpicture}

\thispagestyle{empty} 

\begin{abstract}
This paper focuses on the construction of accurate and predictive data-driven reduced models of large-scale numerical simulations with complex dynamics and sparse training datasets.
In these settings, standard, single-domain approaches may be too inaccurate or may overfit and hence generalize poorly.
Moreover, processing large-scale datasets typically requires significant memory and computing resources which can render single-domain approaches computationally prohibitive.
To address these challenges, we introduce a domain decomposition formulation into the construction of a data-driven reduced model.
In doing so, the basis functions used in the reduced model approximation become localized in space, which can increase the accuracy of the domain-decomposed approximation of the complex dynamics.
The decomposition furthermore reduces the memory and computing requirements to process the underlying large-scale training dataset.
We demonstrate the effectiveness and scalability of our approach in a large-scale three-dimensional unsteady rotating detonation rocket engine simulation scenario with over $75$ million degrees of freedom and a sparse training dataset.
Our results show that compared to the single-domain approach, the domain-decomposed version reduces both the training and prediction errors for pressure by up to $13 \%$ and up to $5\%$ for other key quantities, such as temperature, and fuel and oxidizer mass fractions.
Lastly, our approach decreases the memory requirements for processing by almost a factor of four, which in turn reduces the computing requirements as well.
\end{abstract}

\section{Introduction} \label{sec:intro}
\lettrine{O}ne of the main goals of computational science and engineering is to enable predictive simulations that inform design and decisions in engineering, physics, chemistry or medicine.
To this end, the recent advances in algorithms, computer architectures and high-performance computing enable detailed and accurate such simulations at scale. 
However, making quantitative predictions over long time intervals or performing many-query applications such as design optimization, uncertainty quantification, and control remains challenging even on large supercomputers.

In this paper, we focus on the construction of physics-based data-driven reduced-order models (ROMs) of large-scale numerical simulations characterized by complex dynamics and sparse training datasets, with the goal of replacing the computationally expensive model with a cheap but sufficiently accurate ROM.
The complexity of the underlying dynamics and the sparsity of the training datasets render the construction of accurate and predictive data-driven ROMs challenging via standard, global approaches.
An additional challenge is that processing large-scale training datasets usually requires significant memory and computing resources.
In this work, we address these challenges by introducing a domain decomposition formulation into the construction of data-driven ROMs.

Domain decomposition methods~\cite{DJN15, GK92} are prevalent in scientific computing and especially in high-performance computing where they are used to solve large-scale numerical simulations by decomposing them into smaller pieces that fit, for example, into the memory of one node on a compute cluster or supercomputer.
Domain decomposition has also been used in model reduction.
A large number of approaches target so-called decomposable systems, that is, systems comprising a set of well-defined components.
These approaches are typically based on offline-online decompositions.
In the offline phase, individual reduced bases are constructed for each component.
In the online phase, a ROM for the full system is constructed using domain-decomposition approaches that enforce solution compatibility across interfaces.
The majority of such methods are intrusive and have been developed in the context of linear parametric problems.
Prominent examples include the reduced-basis element (RBE) method proposed in Ref.~\cite{MR02} and the static condensation RBE (scRBE) method~\cite{Ef14, HKP14, HKP13, SP16}.
Reduced bases and domain decomposition were used in Ref.~\cite{IQR16} in the context of networks and complex parametrized geometries.
In recent years, component-based approaches were extended to nonlinear problems as well.
For example, Ref.~\cite{HCC21} proposed an algebraically non-overlapping decomposition based on the least-squares Petrov-Galerkin approach and Ref.~\cite{SD23} formulated a method for nonlinear elliptic partial differential equations (PDEs) based on the partition of unity and adaptive enrichment to update the local reduced spaces.
Ref.~\cite{HDM22} proposed an intrusive component-based approach for large-scale systems, where the components are either ROMs or reduced fidelity full-order models.
Furthermore, the predictive capabilities of the components were enhanced via online basis adaptation~\cite{HD23, Pe20}.
Ref.~\cite{DH23} developed a domain-decomposed ROM approach for coupled heterogeneous systems without assuming a-priori knowledge of the coupled model.

Domain decomposition can also be used to split the full domain into subdomains followed by constructing ROMs for each subdomain.
In the context of computational fluid dynamics simulations, for example, Refs.~\cite{Xi19b, Xi19a} formulated domain-decomposed non-intrusive ROM strategies by leveraging Gaussian process regression or radial basis function interpolation to construct a set of local approximation functions (hypersurfaces) for each subdomain.
Furthermore, Ref.~\cite{An10} used domain decomposition together with balanced truncation for a class of PDE-constrained optimization problems governed by linear time-dependent advection-diffusion equations.
This work was extended in Ref.~\cite{AHH11} to the optimization of linear Stokes problems.
Ref.~\cite{CMB19} used a sparse estimation procedure based on local spatial patches suitable for complex flows.
Other approaches besides domain decomposition for reducing complex, nonlinear problems include reduction to spectral submanifolds~\cite{HP16}, including a recently developed data-driven variant~\cite{Ce22}, and model reduction with online basis adaptation and sampling~\cite{Pe20}. 

Despite these advances, constructing accurate, predictive, and scalable ROMs of large-scale systems with complex dynamics and sparse training datasets remains challenging.
For instance, the underlying simulation codes are complex and often proprietary, which makes the implementation of intrusive approaches difficult.
Moreover, scalability is a major challenge due to the size of the datasets that must be generated for training.
Our goal in the present work is to show that integrating domain decomposition into physics-based data-driven reduced modeling allows a flexible and easy-to-implement construction of accurate and scalable ROMs of large-scale systems.

To construct physics-based data-driven ROMs, we use Operator Inference (OpInf)~\cite{PW16, KPW24}, noting that our domain-decomposed approach is also applicable to other data-driven methods.
As a representative example of a large-scale application that is challenging to reduce via standard single-domain approaches, we consider a three-dimensional unsteady rotating-detonation rocket engine (RDRE) combustion chamber with $75,675,600$ degrees of freedom and a sparse training dataset.
The corresponding high-fidelity simulation is computationally very expensive, requiring  about $1$ million CPU-hours on $16,060$ cores on a supercomputer for $1$ ms of simulated physical time.
The large computational cost prevents the utilization of the high-fidelity model for state predictions over longer time intervals or relevant engineering many-query tasks such as design optimization. 
Moreover, the size of the resulting simulation datasets is usually large, typically on the order of a few hundred gigabytes of disk storage.  
Hence, constructing accurate and scalable ROMs for these complex applications is essential in practice.
We note that due to the complexity of the underlying simulation codes, constructing intrusive ROMs for these problems remains challenging.

The remainder of this paper is organized as follows.
Section~\ref{sec:background} introduces the setup for large-scale numerical simulations and summarizes the OpInf approach for constructing physics-based data-driven ROMs.
Section~\ref{sec:novelty} presents the proposed approach for constructing domain-decomposed ROMs from data, which we instantiate for OpInf.
Section~\ref{sec:results} demonstrates the accuracy gains and scalability of the formulated domain-decomposed approach in a large-scale unsteady $3$D RDRE combustion chamber with a sparse training dataset.
We draw conclusions in Section~\ref{sec:conclusions}.

\section{Single-domain physics-based data-driven reduced modeling} \label{sec:background}
Section~\ref{subsec:SD_FOM} introduces the setup for large-scale high-fidelity simulations.
Section~\ref{subsec:SD_OpInf} summarizes the OpInf approach to learning physics-based ROMs from data.

\subsection{Setup for large-scale high-fidelity simulations} \label{subsec:SD_FOM}
We consider the large-scale system of nonlinear equations defined on the time domain $[t_{\textrm{init}}, t_{\textrm{final}}]$, with $t_{\textrm{init}}$ denoting the initial time and $t_{\textrm{final}}$ the final time,
\begin{equation} \label{eq:FOM_general}
    \dot{\mathbf{q}} = \mathbf{f}(t, \mathbf{q}), \quad \mathbf{q}(t_{\textrm{init}})=\mathbf{q}_{\mathrm{init}},
\end{equation}
where $\mathbf{q}(t) \in \mathbb{R}^{n}$ is the $n$-dimensional vector of state variables at time $t \in [t_{\textrm{init}},  t_{\textrm{final}}]$, $\mathbf{q}_{\mathrm{init}}$ is a specified initial condition, and $\mathbf{f} : [t_{\textrm{init}}, t_{\textrm{final}}] \times \mathbb{R}^{n} \rightarrow \mathbb{R}^{n}$ is a nonlinear function that defines the time evolution of the system state.
The governing equations \eqref{eq:FOM_general} are written as an $n$-dimensional system of ordinary differential equations (ODEs) representing the large-scale system of discretized PDEs over the computational domain denoted by $\Omega \subset \mathbb{R}^d$, where typically $d = 2, 3$.
The state dimension, $n \in \mathbb{N}$, scales with the (large) dimension of the spatial discretization, which we denote by $n_x \in \mathbb{N}$.
Letting $n_s \geq 1$ denote the number of state variables, we have that $n = n_s \times n_x$.

We consider the setup where training data are available over the time interval $[t_\textrm{init}, t_\textrm{train}]$, where by $t_\textrm{train} < t_{\textrm{final}}$ we denote the end of the training time horizon.
We collect $n_t$ snapshots over the training horizon by solving the high-fidelity model \eqref{eq:FOM_general} and recording the high-dimensional state solution at $n_t$ time instants, where we typically have that $n \gg n_t$.
The state solution at time $t_k$ is referred to as the $k$th snapshot and is denoted $\mathbf{q}_k$.
The $n_t$ snapshots are collected into a matrix $\mathbf{Q} \in \mathbb{R}^{n \times n_t}$ with $\mathbf{q}_k$ as its $k$th column:
\begin{equation*}
    \mathbf{Q} =
     \begin{bmatrix}
\vert & \vert & & \vert\\
     \mathbf{q}_1 &
     \mathbf{q}_2 &
     \ldots &
     \mathbf{q}_{n_t}\\
     \vert & \vert & & \vert
     \end{bmatrix}.
\end{equation*}
Given $\mathbf{Q}$, the goal is to construct accurate and predictive ROMs.
In the following, we summarize physics-based data-driven reduced modeling using the OpInf approach.

\subsection{Learning physics-based data-driven reduced models via Operator Inference} \label{subsec:SD_OpInf}
OpInf is a scientific machine learning approach for learning ROMs with polynomial structure from data, where the structure is specified by the governing equations~\cite{PW16, KPW24}. 
For more general types of nonlinearities, lifting transformations can be used to expose (sometimes approximate) polynomial structure in the lifted governing equations~\cite{Qi20}.

The first step in OpInf is to compute the (global) representation of the available snapshots in a low-dimensional subspace. 
This is achieved by first computing the thin singular value decomposition (SVD) of the snapshot matrix $\mathbf{Q}$\footnote{Note that it is often needed to shift and scale $\mathbf{Q}$ prior to computing the global representation of the available snapshots in a low-dimensional subspace, especially in problems with multiple state variables with different scales such as pressure and species mass fractions in reactive flows.},
\begin{equation*}
\mathbf{Q} = \mathbf{V} \mathbf{\Sigma} \mathbf{W}^\top,
\end{equation*}
where $\mathbf{V} \in \mathbb{R}^{n \times n_t}$ contains the left singular vectors, $\mathbf{\Sigma} \in \mathbb{R}^{n_t \times n_t}$ is a diagonal matrix containing the singular values of $\mathbf{Q}$ in non-decreasing order $\sigma_1 \geq \sigma_2 \geq \ldots \geq \sigma_{n_t}$, where $\sigma_j$ denotes the $j$th singular value, and $\mathbf{W} \in \mathbb{R}^{n_t \times n_t}$ contains the right singular vectors.
The first $r \ll n$ columns of $\mathbf{V}$, i.e., the left singular vectors corresponding to the $r$ largest singular values form the POD basis $\mathbf{V}_r \in \mathbb{R}^{n \times r}$.
The low-dimensional representation of the snapshots in the reduced-order linear subspace spanned by $\mathbf{V}_r$ is then obtained by computing
\begin{equation} \label{eq:project_snapshots}
    \hat{\mathbf{Q}} = \mathbf{V}_r^\top \mathbf{Q} \in \mathbb{R}^{r \times n_t}.
\end{equation}

The next step is to determine the reduced operators.
For example, for a ROM with a quadratic state dependency
\begin{equation} \label{eq:ROM_quad_time_cont}
    \dot{\hat{\mathbf{q}}} = \hat{\mathbf{A}}\hat{\mathbf{q}} + \hat{\mathbf{H}}\left(\hat{\mathbf{q}} \otimes \hat{\mathbf{q}} \right),
\end{equation}
we must determine the linear and quadratic reduced operators $\hat{\mathbf{A}} \in \mathbb{R}^{r\times r}$ and $\hat{\mathbf{H}} \in \mathbb{R}^{r\times r^2}$.
OpInf determines the reduced operators that best match the projected snapshot data in a minimum residual sense by solving the following linear least-squares minimization problem
\begin{equation} \label{eq:SD_OpInf_regularization_time_cont}
    \argmin_{\hat{\mathbf{A}}, \hat{\mathbf{H}}} \left\lVert \hat{\mathbf{Q}}^\top\hat{\mathbf{A}}^{\top} + \left(\hat{\mathbf{Q}} \otimes \hat{\mathbf{Q}}\right)^\top \hat{\mathbf{H}}^\top - \dot{\hat{\mathbf{Q}}}^\top \right\rVert_F^2 + \lambda_{\ell} \left\lVert \hat{\mathbf{A}}\right\rVert_F^2 + \lambda_q \left\lVert \hat{\mathbf{H}}\right\rVert_F^2,
\end{equation}
where $F$ denotes the Frobenius norm and $\dot{\hat{\mathbf{Q}}} \in \mathbb{R}^{r \times n_t}$ denotes the time derivative of the reduced snapshot data.
If the high-dimensional time derivatives $\dot{\mathbf{Q}} \in \mathbb{R}^{n \times n_t}$ are provided by the simulation code, $\dot{\hat{\mathbf{Q}}}$ is obtained analogously to~\eqref{eq:project_snapshots} as $\dot{\hat{\mathbf{Q}}}=\mathbf{V}_r^\top \dot{{\mathbf{Q}}}$.
Otherwise, it must be estimated numerically.
In the context of OpInf ROMs for single-element rocket combustors, the works in~\cite{Sc10, QFW21} showed that using higher-order finite difference schemes yields accurate approximations when the training snapshots were not subjected to downsampling.
When downsampling the training snapshots, finite difference approximations may become inaccurate and alternative methods such as time-spectral methods~\cite{SL22} can be utilized instead. 
However, in cases where the snapshots are significantly downsampled, as seen in the RDRE scenario considered in the present paper, accurately approximating $\dot{\hat{\mathbf{Q}}}$ becomes challenging.
Given the sensitivity of least-squares minimizations to perturbations in the right-hand side~\cite{GvL13}, inaccuracies in approximating $\dot{\hat{\mathbf{Q}}}$ can lead to inaccurate ROM predictions.
In such cases, we use the discrete-time version of~\eqref{eq:ROM_quad_time_cont} and a corresponding fully discrete version of~\eqref{eq:SD_OpInf_regularization_time_cont} to derive the reduced operators of the discrete-time ROM. 
In this approach, snapshot time derivatives are no longer needed, and are replaced by time-shifted snapshot data, akin to the dynamic mode decomposition (DMD) technique for learning linear discrete-time systems~\cite{Ku16, Sc10, Tu14}; see also~\cite{Fa23, FMW22}.

Following Ref.~\cite{MHW21}, Tikhonov regularization is added to~\eqref{eq:SD_OpInf_regularization_time_cont} to reduce overfitting and to account for model misspecification or other sources of error, with scalar regularization hyperparameters $\lambda_{\ell}, \lambda_q \in \mathbb{R}$.
The procedure proposed in~\cite{MHW21} for finding the optimal values for $\lambda_{\ell}, \lambda_q$ involves solving~\eqref{eq:SD_OpInf_regularization_time_cont} and ensuing the resulting ROM over a time horizon $[t_{\textrm{init}}, t_{\textrm{reg}}]$ with $t_{\textrm{reg}} \geq t_{\textrm{final}}$ over two nested loops: an outer-loop over candidate values for $\lambda_{\ell}$ and an inner-loop over candidate values for $\lambda_{q}$; see also~\cite{QFW21}.
Nonetheless, this search is embarrassingly parallel and computationally cheap in general since~\eqref{eq:SD_OpInf_regularization_time_cont} depends on the reduced dimension $r \ll n$ and can therefore be solved efficiently via standard least squares methods~\cite{GvL13}.
We note that especially in problems with complex dynamics such as the RDRE scenario considered in this work regularization is key to obtaining accurate OpInf ROMs.

We note a few important points about our OpInf formulation.
Firstly, note that the quadratic reduced operator $\hat{\mathbf{H}} \in \mathbb{R}^{r \times r^2}$ in~\eqref{eq:ROM_quad_time_cont} is not uniquely defined as written due to the symmetry of quadratic products (i.e., $q_1q_2=q_2q_1$), leading to redundant coefficients in the OpInf optimization problem statement~\eqref{eq:SD_OpInf_regularization_time_cont} .
To ensure a unique solution, we eliminate the redundant degrees of freedom, which amounts to learning an operator of dimension $r \times r(r + 1)/2$ as described in \cite{KPW24}.
Secondly, we impose a restriction on the maximum reduced dimension of the OpInf ROM so as to ensure that the OpInf learning problem~\eqref{eq:SD_OpInf_regularization_time_cont} is not under-determined. 
This maximum reduced dimension is set by ensuring that the number of operator coefficients to be learned does not exceed the number of training data points available. 
As discussed in Refs.~\cite{PW16} and \cite{KPW24}, the least squares problem \eqref{eq:SD_OpInf_regularization_time_cont} decouples into $r$ independent least squares problems; nonetheless, in scenarios with sparse training datasets, the sparsity of data can impose restrictions on the ROM dimension.
Another approach might be to permit \eqref{eq:SD_OpInf_regularization_time_cont} to be under-determined and instead employ $\ell_1$-based sparsification instead of Tikhonov regularization, similarly to Ref.~\cite{RLR20}, for example.  
Lastly, we note that additional structure can be imposed on the OpInf ROM. For example, one might impose energy conservation by adding constraints to the entries of the reduced quadratic operator for a quadratic ROM~\cite{LB18}, or impose Hamiltonian~\cite{Sh22, GT23} and Lagrangian~\cite{SK24} structure. 

\section{Domain-decomposed physics-based data-driven reduced modeling} \label{sec:novelty}
This section presents our proposed approach, which integrates domain decomposition into physics-based data-driven reduced modeling of large-scale systems with sparse training datasets.
Section~\ref{subsec:DD_FOM} presents the domain-decomposed formulation of the large-scale governing equations.
This formulation is extended to data-driven reduced modeling and instantiated for OpInf in Section~\ref{subsec:DD_OpInf}.

\subsection{Domain-decomposed formulation of the governing equations} \label{subsec:DD_FOM}
We begin by splitting the computational domain $\Omega$ into $k \in \mathbb{N}_{\geq 2}$ subdomains $\Omega_1, \Omega_2, \ldots, \Omega_k \in \mathbb{R}^d$ such that $\Omega = \cup_{i = 1}^k \Omega_i$ and $k \leq n_x$.
In the limit when $k = n_x$, the computational domain is split into $n_x$ subdomains with a single spatial degree of freedom per subdomain.
We denote by $\mathbf{q}^{(i)}(t) \in \mathbb{R}^{n_i}$ the components of the full-domain vector of state variables $\mathbf{q}(t)$ at time $t$ with size $\mathbb{N} \ni n_i < n$ corresponding to the $i$th subdomain.

For each subdomain $\Omega_i$, let $I(i) = \{j: j \in \{1, 2, \ldots, k\} \setminus \{i\}, \Omega_j \text{is  adjacent to } \Omega_i \}$ denote the set of indices of all subdomains adjacent to it.  
We can then rewrite the governing equations \eqref{eq:FOM_general} as a system of coupled ODEs: 
\begin{equation} \label{eq:DD_FOM}
\begin{cases}
   \dot{\mathbf{q}}^{(1)} & = \underbrace{\mathbf{f}_{1, 1}(t, \mathbf{q}^{(1)})}_{\mathrm{interior} \ \Omega_1} + \underbrace{\sum_{j \in I(1)} \mathbf{f}_{1, j}(t, \mathbf{q}^{(j)})}_{\mathrm{coupling} \ \Omega_1 \rightarrow \Omega_j} \\
   \dot{\mathbf{q}}^{(2)} & = \underbrace{\mathbf{f}_{2, 2}(t, \mathbf{q}^{(2)})}_{\mathrm{interior} \ \Omega_2} + \underbrace{\sum_{j \in I(2)} \mathbf{f}_{2, j}(t, \mathbf{q}^{(j)})}_{\mathrm{coupling} \ \Omega_2 \rightarrow \Omega_j} \\
   \vdots \\
    \\
    \dot{\mathbf{q}}^{(k)} & = \underbrace{\mathbf{f}_{k, k}(t, \mathbf{q}^{(k)})}_{\mathrm{interior} \ \Omega_k} + \underbrace{\sum_{j \in I(k)} \mathbf{f}_{k, j}(t, \mathbf{q}^{(j)})}_{\mathrm{coupling} \ \Omega_k \rightarrow \Omega_j},
    \end{cases}
\end{equation}
where $\mathbf{f}_{i, i} : [t_{\textrm{init}}, t_{\textrm{end}}] \times \mathbb{R}^{n_i} \rightarrow \mathbb{R}^{n_i}$ models the time evolution of the state in the interior of the $i$th subdomain and $\mathbf{f}_{i, j}: [t_{\textrm{init}}, t_{\textrm{end}}] \times \mathbb{R}^{n_j} \rightarrow \mathbb{R}^{n_i}$ models the state evolution at the interface between two adjacent subdomains $i$ and $j$.
Note that we did not specify whether the decomposition is done with or without overlapping; we will do so in the next section.

\subsection{Domain-decomposed Operator Inference} \label{subsec:DD_OpInf}
\subsubsection{General considerations}
Starting from the domain-decomposed formulation~\eqref{eq:DD_FOM} of the governing equations, we seek the construction of accurate, predictive, and scalable physics-based domain-decomposed ROMs from data.
The generic steps to construct such ROMs are to (i) split the training dataset into subsets corresponding to subdomains, (ii) use these subsets to construct reduced bases for each subdomain, and (iii) determine the corresponding interior and coupling reduced operators. 
Once constructed, the reduced system of coupled ODEs is integrated over the desired time horizon to obtain the reduced solutions for each subdomain, which are then used to determine the approximate global, full-domain solution.
Constructing accurate data-driven domain-decomposed ROMs using a non-overlapping decomposition can be challenging in practice, especially in problems characterized by complex dynamics such as the RDRE scenario considered in Section~\ref{sec:results}.
In such problems, it is likely that the corresponding reduced solutions will have sharp transitions or oscillations at the interfaces even when the reduced dimensions of the subdomain ROMs are large, which in turn will lead to inaccurate global solutions (see Refs.~\cite{BCI13, HDM22}, for example, for similar discussions).
To mitigate this challenge, we formulate our domain-decomposed approach in terms of overlapping subdomains.
That is, for two adjacent subdomains $i$ and $j$, we have that $\Omega_i \cap \Omega_j = \Omega_{i \cap j}$ where $\Omega_{i \cap j}$ denotes the overlapping region in the physical space.
Notice that overlapping implies that $\sum_{i = 1}^k n_i > n$.

The decomposition can be influenced by several factors, including the domain topology, the anisotropy of the dynamics, and the size of the training dataset.
For instance, if the computational domain consists of multiple similar components, it is natural to decompose the domain into subdomains corresponding to each component. 
An illustration of this approach can be found in Ref.~\cite{HDM22}, where a multi-injector rocket combustor was decomposed with each subdomain representing an injector element.
Additionally, if the complex dynamics are concentrated in specific regions of the domain, it can be advantageous to incorporate this localization into the decomposition process. This approach was considered in Ref.~\cite{LKB03}, for example, where the computational domain was decomposed based on the shock location. 
Another crucial consideration is the size of the training dataset. 
When dealing with a large training dataset, it is essential to consider the available computational resources for managing the dataset effectively, which in turn may dictate a required level of partitioning.

In the following, we present our formulation for data-driven domain-decomposed ROMs with overlapping for OpInf, noting that it can be straightforwardly applied to other data-driven approaches such as DMD. 
We use the notation \emph{DD-$k$} to refer to the decomposition of $\Omega$ into $k$ overlapping subdomains and \emph{DD-$k$ OpInf} to refer to the corresponding domain-decomposed OpInf formulation. 

\subsubsection{Learning the domain-decomposed reduced operators via Operator Inference} \label{subsubsec:DD_OpInf_steps}
Based on the overlapping decomposition of $\Omega$ into $\Omega_1, \Omega_2, \ldots, \Omega_k$, we split the snapshot matrix $\mathbf{Q}$ into $k$ snapshot matrices $\mathbf{Q}^{(i)} \in \mathbb{R}^{n_i \times n_t}, \quad i = 1, 2, \ldots, k$.
In the first step in DD-$k$ OpInf, we compute the representation of these snapshots in low-dimensional subspaces for all $k$ subdomains\footnote{As noted in Section~\ref{subsec:SD_OpInf}, especially in problems with multiple state variables with different scales, we need to first shift and scale the given snapshot data prior to performing the overlapping decomposition and computing the POD basis in each subdomain.}.
This is achieved by first computing the thin SVD of $\mathbf{Q}^{(i)}$,
\begin{equation*}
\mathbf{Q}^{(i)} = \mathbf{V}^{(i)} \mathbf{\Sigma}^{(i)} \left(\mathbf{W}^{(i)}\right)^\top,
\end{equation*}
where $\mathbf{V}^{(i)} \in \mathbb{R}^{n_i \times n_t}$, $\mathbf{\Sigma}^{(i)} \in \mathbb{R}^{n_t \times n_t}$, and $\mathbf{W}^{(i)} \in \mathbb{R}^{n_t \times n_t}$.
The left singular vectors corresponding to the $r_i \ll n_i$ largest singular values form the POD basis $\mathbf{V}^{(i)}_{r_i} \in \mathbb{R}^{n_i \times r_i}$ for the $i$th subdomain.
The low-dimensional representation of the snapshots in the reduced-order subspace spanned by $\mathbf{V}^{(i)}_{r_i}$ is then obtained by computing
\begin{equation} \label{eq:DD_OpInf_projections}
    \hat{\mathbf{Q}}^{(i)} = \left(\mathbf{V}^{(i)}_{r_i}\right)^\top \mathbf{Q}^{(i)} \in \mathbb{R}^{r_i \times n_t}.
\end{equation}
For the continuous-time OpInf formulation, we must also estimate $\dot{\hat{\mathbf{Q}}}^{(i)} \in \mathbb{R}^{r_i \times n_t}$.
\begin{remark}
The size of the overlapping region has a direct effect on the corresponding POD bases in that it influences how well these bases capture the dynamics in the neighboring subdomains.
For two adjacent subdomains with similar dynamics, a small overlapping region generally suffices.
In contrast, a larger overlapping region is typically needed for neighboring subdomains with anisotropic dynamics to capture potentially sharp transitions, for example.
\end{remark}

We note that if the dynamics within the individual subdomains exhibit self-similarity, it may be feasible to compute POD bases for a single domain or a few subdomains and subsequently use these bases for the remaining subdomains.
For example, Ref.~\cite{HDM22} considered such an approach in the context of intrusive ROMs for rocket combustors with multiple injector elements.
In this approach, one would need to accommodate potential dynamics variations from subdomain to subdomain. 
One strategy would be to parametrize the boundaries of the subdomain(s) for which the POD bases are computed;
however, finding such a parametrization for problems with complex dynamics such as the RDRE scenario considered in our numerical experiments in Section~\ref{sec:results} is not straightforward. 
We therefore do not pursue this idea in this paper, but rather we compute a POD basis for each subdomain individually.

For a single-domain ROM with a quadratic state dependency (cf.~\eqref{eq:ROM_quad_time_cont}) and linear coupling terms, the corresponding DD-$k$ OpInf ROM reads
\begin{equation} \label{eq:DD_OpInf_time_cont}
\begin{cases}
   \dot{\hat{\mathbf{q}}}^{(1)} & = \underbrace{\hat{\mathbf{A}}_{1, 1}\hat{\mathbf{q}}^{(1)} + \hat{\mathbf{H}}_{1, 1}\left(\hat{\mathbf{q}}^{(1)} \otimes \hat{\mathbf{q}}^{(1)} \right)}_{\textrm{interior}} + \underbrace{\sum_{j \in I(1)} \hat{\mathbf{A}}_{1, j}\hat{\mathbf{q}}^{(j)}}_{\textrm{coupling}} \\
   \dot{\hat{\mathbf{q}}}^{(2)} & = \underbrace{\hat{\mathbf{A}}_{2, 2}\hat{\mathbf{q}}^{(2)} + \hat{\mathbf{H}}_{2, 2}\left(\hat{\mathbf{q}}^{(2)} \otimes \hat{\mathbf{q}}^{(2)} \right) }_{\textrm{interior}} + \underbrace{\sum_{j \in I(2)} \hat{\mathbf{A}}_{2, j}\hat{\mathbf{q}}^{(j)}}_{\textrm{coupling}} \\
   \vdots \\
    \\
    \dot{\hat{\mathbf{q}}}^{(k)} & = \underbrace{\hat{\mathbf{A}}_{k, k}\hat{\mathbf{q}}^{(k)} + \hat{\mathbf{H}}_{k, k}\left(\hat{\mathbf{q}}^{(k)} \otimes \hat{\mathbf{q}}^{(k)} \right)}_{\textrm{interior}} + \underbrace{\sum_{j \in I(k)} \hat{\mathbf{A}}_{k, j}\hat{\mathbf{q}}^{(j)}}_{\textrm{coupling}}.
    \end{cases}
\end{equation}
The linear coupling terms in~\eqref{eq:DD_OpInf_time_cont} act like forcing/boundary condition terms between neighboring subdomains.
\begin{figure}[htb!]
\centering
\includegraphics[width=0.7\textwidth]{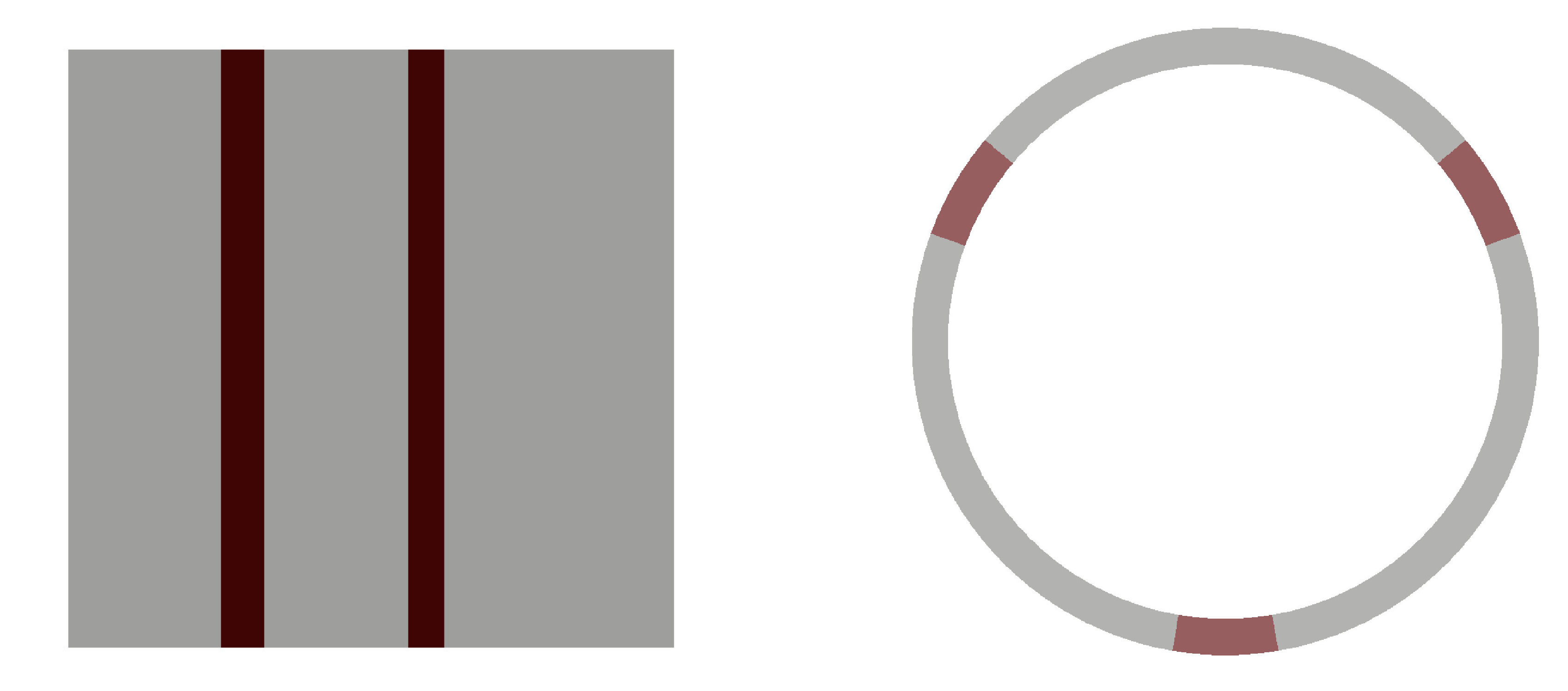}
\caption{Two example domains decomposed into three overlapping subdomains. The overlapping regions are shown in red.}
\label{fig:DD_ex}
\end{figure}
We note that the coupling between adjacent subdomains may also include additional components such as quadratic terms. 
Extending~\eqref{eq:DD_OpInf_time_cont} to such settings is trivial.
Figure~\ref{fig:DD_ex} provides a visual illustration of two example domains.
Both domains are decomposed into three overlapping subdomains.
The left plot shows a rectangular domain in which, from left to right, the first subdomain is coupled with the second subdomain, and the second subdomain is coupled with the third.
Therefore, $I(1) = \{2\}, I(2) = \{1, 3\}$, and $I(3) = \{2\}$.
On the right, we plot an annular, periodic domain in which all three sectors overlap, which means that $I(1) = \{2, 3\}, I(2) = \{1, 3\}$, and $I(3) = \{1, 2\}$.

In the next step, we learn the reduced operators for each subdomain in~\eqref{eq:DD_OpInf_time_cont} via OpInf by solving $k$ independent linear least-squares minimizations.
That is, in the $i$th least-squares minimization problem with $i = 1, 2, \ldots, k$, we learn the interior linear and quadratic reduced operators $\hat{\mathbf{A}}_{i, i} \in \mathbb{R}^{r_i \times r_i}$ and $\hat{\mathbf{H}}_{i, i} \in \mathbb{R}^{r_i \times r_i^2}$ as well as the coupling operators $\hat{\mathbf{A}}_{i, j} \in \mathbb{R}^{r_i \times r_j}$ for all $j \in I(i)$ by solving
\begin{multline} \label{eq:DD_OpInf_subdom_regularization_time_cont}
    \argmin_{\hat{\mathbf{A}}_{i, i}, \hat{\mathbf{H}}_{i, i}, \hat{\mathbf{A}}_{i, j}, \quad j \in I(i)} \left\lVert  \left(\hat{\mathbf{Q}}^{(i)}\right)^\top \hat{\mathbf{A}}^{\top}_{i, i} + \left(\hat{\mathbf{Q}}^{(i)} \otimes \hat{\mathbf{Q}}^{(i)}\right)^\top \hat{\mathbf{H}}^{\top}_{i, i} + \sum_{j \in I(i)}\left(\hat{\mathbf{Q}}^{(j)}\right)^\top \hat{\mathbf{A}}^{\top}_{i, j} - \left(\dot{\hat{\mathbf{Q}}}^{(i)}\right)^\top \right\rVert_F^2 + \\ \lambda_{\ell}^{(i)} \left( \left\lVert \hat{\mathbf{A}}_{i, i}\right\rVert_F^2 + \sum_{j \in I(i)} \left\lVert \hat{\mathbf{A}}_{i, j}\right\rVert_F^2 \right) + \lambda_q^{(i)} \left\lVert \hat{\mathbf{H}}_{i, j}\right\rVert_F^2.
\end{multline}
Notice that we use different regularization hyperparameters $\lambda_{\ell}^{(i)}, \lambda_q^{(i)} \in \mathbb{R}$ for each subdomain.

To find the optimal values for $\{\lambda_{\ell}^{(i)}, \lambda_q^{(i)}\}_{i=1}^k$, we extend the strategy formulated in Ref.~\cite{MHW21} (cf.~Section~\ref{subsec:SD_OpInf}) to the domain-decomposed setting.
Since the reduced system of ODEs~\eqref{eq:DD_OpInf_time_cont} is coupled, this involves solving~\eqref{eq:DD_OpInf_time_cont} over $2k$ nested loops, that is, two loops for each subdomain for a ROM with quadratic structure. 
Therefore, the effort for finding separate hyperparameters for each subdomain scales exponentially with $2k$ for $k$ overlapping subdomains.
Even though this search is embarrassingly parallel, the overall computational cost can become prohibitive for a large number of subdomains.
An alternative approach in such cases is to instead find a single pair of hyperparameters for all subdomains, keeping in mind that this strategy may impact the prediction accuracy of the ensuing DD-$k$ OpInf ROM, especially when the dynamics in the underlying subdomains are heterogeneous.
Once the values of the regularization hyperparameters are determined, the least-square minimization problems~\eqref{eq:DD_OpInf_subdom_regularization_time_cont} are solved independently to find the interior and coupling reduced operators for all $k$ subdomains.

In summary, the main ingredients needed in our DD-$k$ OpInf approach are (i)~knowledge of the structure of the governing equations, (ii)~simulation data from the high-fidelity simulation code, and (iii)~the coordinates of the underlying computational domain to perform the overlapping decomposition.
From an implementation perspective, the proposed approach is flexible and easy to implement.
After the decomposition is performed, the major required operations are the thin SVD (or equivalent approaches) of the subdomain snapshot matrices to determine the POD bases and linear least-squares minimization to learn the reduced interior and coupling operators.
These are standard numerical linear algebra operations with efficient implementations provided by most scientific computing libraries.

\subsubsection{Obtaining the full-domain solution using the domain-decomposed reduced solutions} \label{subsubsec:DD_recon}
Once the interior and coupling reduced operators are learned, the next step is to integrate the reduced system of coupled ODEs~\eqref{eq:DD_OpInf_time_cont} over the time horizon $[t_{\textrm{init}}, t_{\textrm{final}}]$ and collect the $k$ reduced solutions for all subdomains.
The reduced initial conditions $\hat{\mathbf{q}}^{(i)}_{\textrm{init}} \in \mathbb{R}^{r_i}$ for each subdomain are determined by projecting their corresponding full-domain initial conditions $\mathbf{q}^{(i)}_{\textrm{init}} \in \mathbb{R}^{n_i}$ onto the linear subspace spanned by the column vectors of $\mathbf{V}^{(i)}_{r_i}$ for all $i = 1, 2, \ldots, k$:
\begin{equation*}
    \hat{\mathbf{q}}^{(i)}_{\textrm{init}} = \left(\mathbf{V}^{(i)}_{r_i}\right)^\top \mathbf{q}^{(i)}_{\textrm{init}}.
\end{equation*}
We denote by $\tilde{\mathbf{q}}^{(i)}(t) \in \mathbb{R}^{r_i}$ the reduced solution at time $t$ corresponding to the $i$th subdomain.
We then map $\tilde{\mathbf{q}}^{(i)}(t)$ back to the original space to obtain the approximate DD-$k$ OpInf solutions in the original coordinates $\mathbf{q}^{(i)}_{\textrm{DD}-k}(t) \in \mathbb{R}^{n_i}$: 
\begin{equation*}
    \mathbf{q}^{(i)}_{\textrm{DD}-k}(t) = \mathbf{V}^{(i)}_{r_i} \tilde{\mathbf{q}}^{(i)}(t) \in \mathbb{R}^{n_i}, \quad i = 1, 2, \ldots, k.
\end{equation*}

To determine the full-domain approximate solution $\mathbf{q}_{\textrm{DD}-k}(t) \in \mathbb{R}^{n}$ using $\mathbf{q}^{(1)}_{\textrm{DD}-k}(t), \mathbf{q}^{(2)}_{\textrm{DD}-k}(t), \ldots, \mathbf{q}^{(k)}_{\textrm{DD}-k}(t)$, we proceed as follows.
The interior degrees of freedom in each subdomain solution are uniquely mapped to their corresponding indices in the full-domain solution.
The overlapping degrees of freedom, however, cannot be uniquely mapped since they are computed twice for two adjacent subdomains.
To obtain the full-domain solution, we smoothly combine the overlapping degrees of freedom.
We illustrate this process in the following example.

\begin{figure}[htb!]
\centering
\includegraphics[width=1.0\textwidth]{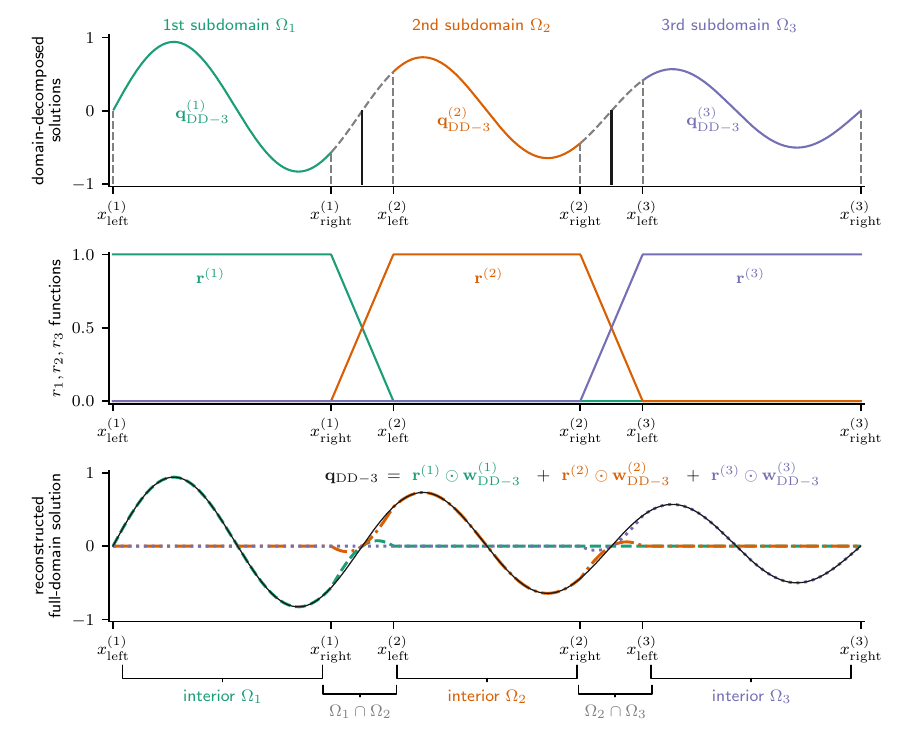}
\caption{Visual illustration of the process of determining the full-domain solution using domain-decomposed solutions using an example with three overlapping subdomains.}
\label{fig:DD_recon_ex}
\end{figure}

Consider a stationary damped sinusoidal wave function defined on a one-dimensional spatial domain $\Omega \subset \mathbb{R}$.
We split $\Omega$ into three overlapping subdomains $\Omega_1, \Omega_2$, and $\Omega_3$ such that $\Omega_1$ overlaps with $\Omega_2$ and the $\Omega_2$ overlaps with $\Omega_3$.
This is illustrated in the top part in Figure~\ref{fig:DD_recon_ex}.
The subdomains are separated by the black solid vertical lines.
Moreover, let $\mathbf{q} \in \mathbb{R}^n$ denote the function evaluated at $n$ spatial points and let $\mathbf{q}^{(i)}_{\textrm{DD}-3} \in \mathbb{R}^{n_i}$ denote the respective subdomain evaluations for $i = 1, 2, 3$.
In addition, we denote by $\{ x^{(i)}_{\textrm{left}}, x^{(i)}_{\textrm{right}} \}_{i=1}^3$ the spatial locations delimiting the interior of the three subdomains (shown on the $x$-axes).

To obtain the full-domain solution $\mathbf{q}_{\textrm{DD}-3} \in \mathbb{R}^{n}$ using $\mathbf{q}^{(1)}_{\textrm{DD}-3}$, $\mathbf{q}^{(2)}_{\textrm{DD}-3}$, and $\mathbf{q}^{(3)}_{\textrm{DD}-3}$, we first construct three functions $r^{(i)} : \Omega \rightarrow [0, 1]$ -- one for for each subdomain -- that keep the respective interior degrees of freedom unchanged and at the same time permit us to smoothly combine the corresponding overlapping degrees of freedom:
\begin{equation} \label{eq:DD_recon_func}
r^{(i)}(x) = 
\begin{cases}
    1, & x^{(i)}_{\textrm{left}} \leq x \leq x^{(i)}_{\textrm{right}} \\
    \left(x - x^{(i + 1)}_{\textrm{left}}\right)/\left(x^{(i)}_{\textrm{right}} - x^{(i + 1)}_{\textrm{left}}\right), & x^{(i)}_{\textrm{right}} < x < x^{(i + 1)}_{\textrm{left}}, \quad i = 1, 2  \\
    \left(x - x^{(i - 1)}_{\textrm{right}}\right)/\left(x^{(i)}_{\textrm{left}} - x^{(i - 1)}_{\textrm{right}}\right), & x^{(i)}_{\textrm{left}} < x < x^{(i - 1)}_{\textrm{right}}, \quad i = 2, 3  \\
    0, & \textrm{otherwise}.
\end{cases}
\end{equation}
We then evaluate these functions at all spatial points in the full-domain to respectively obtain $\mathbf{r}^{(1)}, \mathbf{r}^{(2)}$, and $\mathbf{r}^{(3)} \in \mathbb{R}^{n}$, which are depicted in the center plot in Figure~\ref{fig:DD_recon_ex}.
We extend $\mathbf{q}^{(1)}_{\textrm{DD}-3}$, $\mathbf{q}^{(2)}_{\textrm{DD}-3}$, and $\mathbf{q}^{(3)}_{\textrm{DD}-3}$ to the full-domain (by adding zeros to all additional entries) to obtain $\mathbf{w}^{(1)}_{\textrm{DD}-3}, \mathbf{w}^{(2)}_{\textrm{DD}-3}, \mathbf{w}^{(3)}_{\textrm{DD}-3} \in \mathbb{R}^n$. 
The full-domain solution $\mathbf{q}_{\textrm{DD}-3}$ is computed as:
\begin{equation} \label{eq:DD_full_dom_rec}
    \mathbf{q}_{\textrm{DD}-3} = \sum_{i = 1}^{3} \mathbf{r}^{(i)} \odot \mathbf{w}^{(i)}_{\textrm{DD}-3},
\end{equation}
where $\odot$ denotes the element-wise (Hadamard) product.
This is depicted in the bottom plot in Figure~\ref{fig:DD_recon_ex}. 
Computing the full-domain solution requires knowledge of the solution in each subdomain, but it does not require those subdomain solutions to take any specific functional form.
It is trivial to extend this procedure to an arbitrary number of subdomains as well as to other types of domains such as periodic domains by suitably adjusting the definition of the functions in~\eqref{eq:DD_recon_func}.

\section{Numerical demonstration in a large-scale rotating detonation rocket engine simulation} \label{sec:results}
We demonstrate the capabilities of the presented domain-decomposed ROM approach in a large-scale three-dimensional unsteady RDRE scenario with $75,675,600$ degrees of freedom and a sparse training dataset.
RDREs are based on fuel injection in an axially symmetric chamber such as an annulus wherein, once ignited at suitable conditions, a system of spinning detonation waves is produced (see Ref.~\cite{By06} for a detailed description).
This concept, in addition to its mechanical simplicity, offers several other propulsive advantages compared to more conventional rocket engine designs, making RDREs an active area of research in recent years.
However, the large computational cost of the corresponding simulations prohibits using the high-fidelity model for state predictions over long time intervals or for practically relevant engineering applications such as the design optimization of these devices.
This motivates the development of accurate and predictive data-driven ROMs to ultimately enable these tasks at scale.

\subsection{Problem setup for high-fidelity large-eddy simulations}
The physics of the problem are modeled using the three-dimensional, reactive, viscous Navier-Stokes equations coupled with a skeletal chemistry-mechanism (FFCMy-12) based on the FFCM model~\cite{FFCM1}.
The full RDRE simulation was performed using implicit large-eddy simulations (LES) via the AHFM (ALREST High-Fidelity Modeling) simulation code from the Air Force Research Laboratory (AFRL); for more details, please refer to Refs.~\cite{Ba20, Li19}.
The computational mesh size for the full RDRE comprises around $136$ million spatial cells.
The time step, computed adaptively, is about one nanosecond in the quasi-steady state regime.
The corresponding high-fidelity LES is computationally expensive, requiring about $1$ million CPU-hours on $16,060$ cores on a supercomputer to obtain $1$ ms of simulated physical time.
Due to the size of the resulting simulation datasets, the number of time instants for which the high-fidelity LES solutions can be saved to disk is limited, usually amounting to only a few hundred.

The scenario considered here comprises $72$ discrete injector-pairs.
The values of the input parameters characterizing the flow conditions, the mass flow-rate $\dot{\textrm{m}} \ [\textrm{kg} \cdot \textrm{s}^{-1}]$ and equivalence ratio $\Phi$, are $\dot{\textrm{m}} = 0.267 \ \textrm{kg} \cdot \textrm{s}^{-1}$ and $\Phi = 1.16$.
These values lead to the formation of two dominant co-rotating waves and no secondary waves in the quasi-limit cycle regime.
We note that this setup, utilizing a sparser dataset, was also considered in Ref.~\cite{Fa23} within the framework of data-driven parametric OpInf ROMs for RDRE combustion chambers.

\subsection{Setup for data-driven reduced modeling}
\begin{figure}[htb!]
\centering
\includegraphics[width=0.8\textwidth]{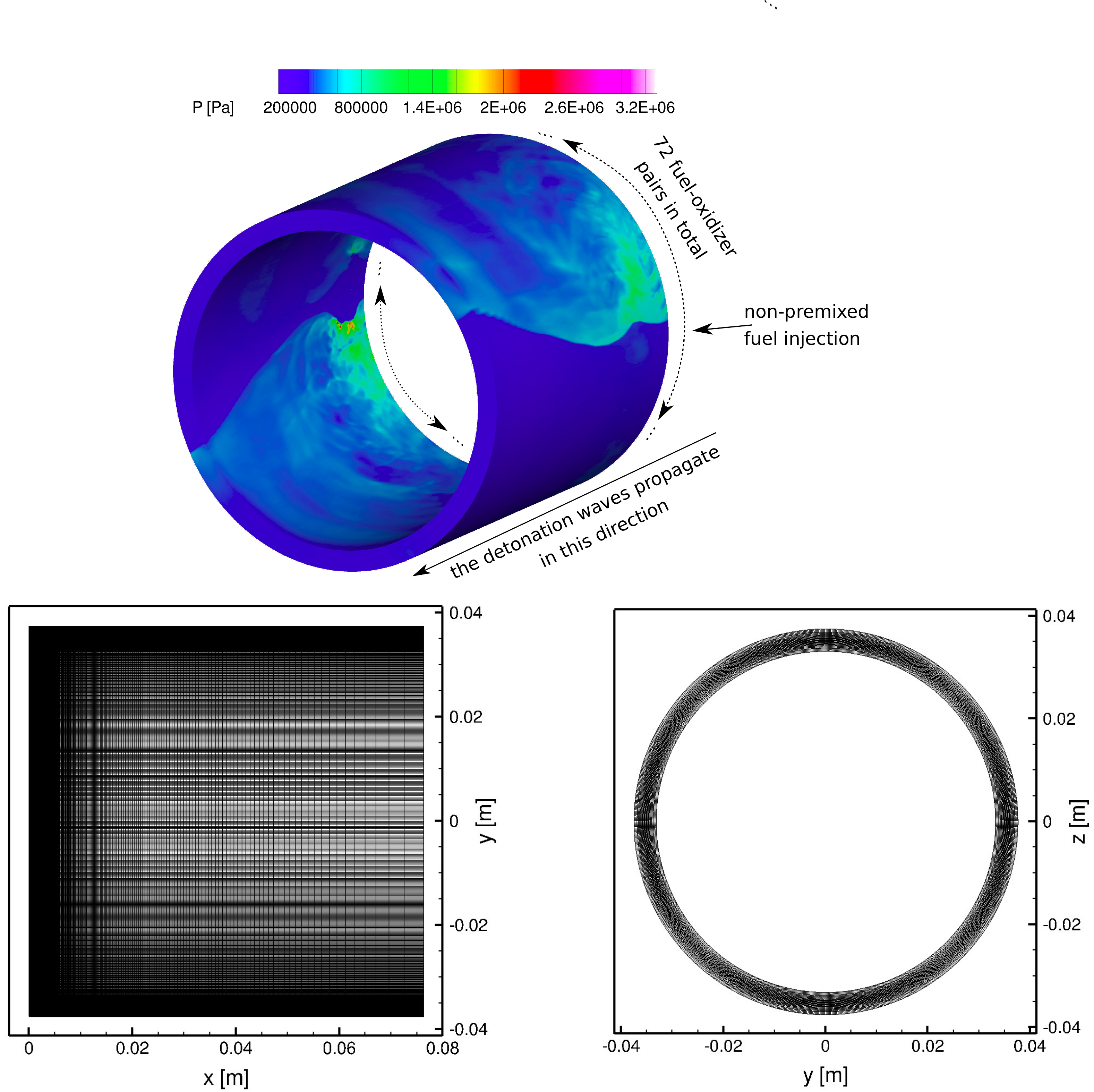}
\caption{Pressure field example (top), and the $x$-$y$ and $y$-$z$ extents of the computational domain (bottom). For more details, we refer the reader to Ref.~\cite{Be23}.}
\label{fig:RDE_domain_and_pressure_ex}
\end{figure}

We are interested in constructing ROMs for the combustion chamber only.
The ROM domain spans from $0.05$ to $76.15$ mm in the $x$ direction, spans from $-37.50$ to $37.50$ mm in both $y$ and $z$ directions, and has a fixed channel height of $4.44$ mm throughout.
The original simulation data has been interpolated onto a grid domain comprising $n_x = 4,204,200$ spatial degrees of freedom, with clustering of grid points at mid channel and closer toward the injector plane.
Figure~\ref{fig:RDE_domain_and_pressure_ex} plots a typical pressure field example (top) as well as the $x$-$y$ and $y$-$z$ extents of the domain (bottom).
For more details about this scenario, including the geometry, we refer the reader to Ref.~\cite{Be23}.

Following the work in Ref.~\cite{Sw20}, we build the ROMs for the RDRE combustion chamber using the specific volume flow variables (specific volume, pressure and the three velocity components), chemical species mass fractions, and temperature, which makes most terms in the governing equations linear or quadratic. 
The data snapshots are thus transformed to represent the following $n_s = 18$ state variables:    
\begin{equation} \label{eq:RDE_state_vars}
    q = [1/\rho \ \ p \ \  v_x \ \ v_y \ \ v_z \ \ w_{\mathrm{CH_4}} \ \ w_{\mathrm{O_2}} \ \ w_{\mathrm{CO_2}} \ \ w_{\mathrm{H_2O}} \ \ w_{\mathrm{CO}} \ \ w_{\mathrm{H_2}} \ \ w_{\mathrm{OH}} \ \ w_{\mathrm{CH_2O}} \ \ w_{\mathrm{CH_4}} \ \ w_{\mathrm{HO_2}} \ \ w_{\mathrm{H}} \ \ w_{\mathrm{O}} \ \ T]^\top,
\end{equation}
where $\rho \ [\mathrm{kg} \cdot \mathrm{m}^{-3}]$ is density (and $1/\rho \ [\mathrm{m}^{3} \cdot \mathrm{kg}^{-1}]$ is specific volume), $p \ [\mathrm{Pa}]$ is pressure, $v_x, v_y, v_z \ [\mathrm{m} \cdot \mathrm{s}^{-1}]$ are the velocity components, $w_{\cdot}$ are the mass fractions of all chemical species ($12$ in total here), and $T \ [\mathrm{K}]$ is temperature.
The total number of degrees of freedom is therefore $n = n_s \times n_x = 18 \times 4,204,200 = 75,675,600$.
Since the $n_s = 18$ state variables in~\eqref{eq:RDE_state_vars} have significantly different scales (for example, from the order of $10^4$ to $10^6$ Pa for pressure to between $0$ and $1$ for the species mass fractions), we center and scale the snapshot data variable-by-variable prior to constructing ROMs.
Each state variable is first centered around the mean field (over the training data) in that variable, and then scaled by its maximum absolute value, thus ensuring that the scaled variables do not exceed $[-1, 1]$.

\begin{remark}
Centering and scaling of snapshots are essential for obtaining accurate OpInf ROMs, in order to prevent undesirable bias towards variables with larger values in the corresponding POD bases.
Scaling is particularly critical in cases where the state comprises multiple physical quantities with differing physical scales.
In this work, we perform centering around the mean value of each state variable over the training horizon.
We explored two scaling parameters: maximum absolute value and standard deviation of centered variables. 
Our findings show that, for this particular problem, the choice of scaling parameters from these two options does not have a significant impact on OpInf ROM prediction quality. 
The key consideration is to maintain consistent scales for non-dimensionalized variables.
\end{remark}

\begin{remark}
In some applications the choice of scaling parameters can have a large effect, and for our application there may be a better scaling choice that would lead to improved ROM performance.
It remains an open and important challenge to identify optimal scaling strategies. A particularly interesting avenue of work would be to explore strategies that non-dimensionalize the snapshots in a physically consistent way. 
Non-dimensionalization and variable transformations are well known to play a key role in establishing desirable numerical properties in high-fidelity computational fluid dynamics simulations; their role in reduced-order modeling is much less explored.
As an alternative to an explicit variable scaling, previous work has explored physically-motivated choices for inner products in reduced-order modeling, such as an energy-based inner product for compressible flows in Ref.~\cite{rowley2004model} and an entropy normalization scheme that is related to the Chu energy norm in Ref.~\cite{vogel2022novel}.
For further discussions on data normalization beyond non-intrusive data-driven ROMs, we refer to Refs.~\cite{HDM22, LP97}.
\end{remark}

We have $501$ down-sampled snapshots in $\mathbb{R}^{75,675,600}$ from the high-fidelity LES over the time interval $[2.50, 3.75]$ ms, which corresponds to roughly four full-cycle periods of the two-wave system. 
The time step between the down-sampled snapshots is about $2.5$ $\upmu$s, which means that every $2500$th time instant or so from the high-fidelity simulation was saved to disk.
In our numerical experiments, we use the first $n_t = 375$ snapshots to train the ROMs, which correspond to the about the first three full cycles.
The remaining fourth full cycle, comprising $126$ snapshots, is used to assess the ROMs accuracy for predictions beyond the training horizon.
 
\subsection{Results using domain-decomposed Operator Inference} \label{subsec:RDE_results}
In the following, we use the DD-$k$ OpInf approach presented in Section~\ref{subsec:DD_OpInf} to construct domain-decomposed ROMs and compare it in terms of scalability and accuracy with the standard OpInf method applied to a single domain as summarized in Section~\ref{subsec:SD_OpInf}. 
Here, we refer to single-domain OpInf as \emph{SD OpInf}.
We note that due to the significant downsampling of the available snapshots, it is difficult to accurately estimate the time derivatives of the projected snapshots in both SD and DD-$k$ OpInf.
We therefore instead use the fully discrete version of OpInf where we learn a fully discrete ROM, rather than a semi-discrete ROM; see also~\cite{Fa23, FMW22}.
All ROM calculations were performed in standard double precision arithmetic using a shared-memory machine with $256$ AMD EPYC 7702 CPUs and $2$ TB of RAM, using the \texttt{numpy} and \texttt{scipy} scientific computing libraries in \texttt{python}.

\begin{figure}[htb!]
\centering
\includegraphics[width=0.4\textwidth]{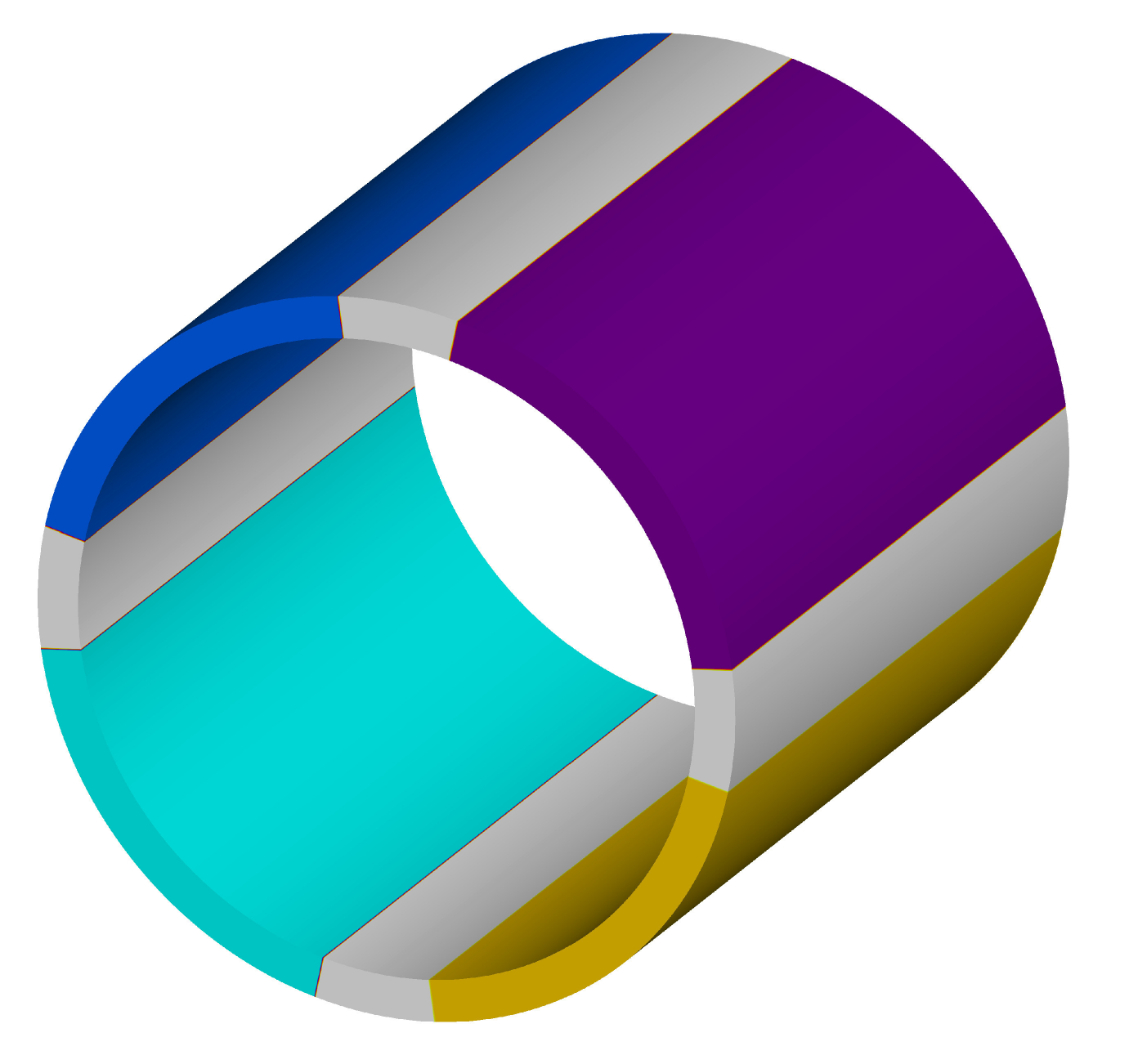}
\caption{Combustion chamber decomposition into four overlapping sector with overlapping region between two sectors of $\pi/9$ radians.}
\label{fig:RDE_DD_split}
\end{figure}

We decompose the $360^{\circ}$ annular combustion chamber into $k=4$ overlapping sectors.
Due to the periodicity of the underlying annular domain, we overlap the first and fourth sectors as well.
This is plotted in Figure~\ref{fig:RDE_DD_split}.
Our experiments using two angle extents for the overlapping region, $2\pi/36 = \pi/18$ and $2\pi/18 = \pi/9$, indicated that the larger overlapping region, $\pi/9$, corresponding to the angle extent of four injectors, leads to the most accurate results.
Therefore, in the following, we consider an overlapping region spanning $\pi/9$ radians. 
\begin{remark}
We have also performed numerical experiments in which we split the domain into $k = 8, 12$, and $24$ overlapping sectors, using several values for the angle extent of the overlapping region as well.
However, we did not obtain any gains in accuracy compared to the case with $k = 4$ overlapping sectors.
This is because for these splits, we used a single, global set of hyperparameters for all $k$ least-squares minimization problems~\eqref{eq:DD_OpInf_subdom_regularization_time_cont} since it was computationally impractical to find separate regularization hyperparameters (see the discussion in Section~\ref{subsubsec:DD_OpInf_steps}).
Since we observed that the DD-$k$ OpInf ROM accuracy depends sensitively on the choice of regularization parameters, using the same regularization hyperparameters rather than individual ones had a detrimental effect on the ROM accuracy.
We therefore consider only a decomposition with $k=4$ overlapping sectors.
We note, nevertheless, that for different configurations, flow conditions, or number of detonation waves, decomposing the annular combustion chamber into a larger number of sectors could increase the accuracy of the domain-decomposed ROM.
\end{remark}

In general, for an RDRE scenario with $N$ dominant waves, it is natural to decompose the annular combustion chamber into $k = 2^{\eta} N$ overlapping sectors, where $\eta \in \mathbb{N}_0$.
Moreover, for a configuration with $M$ injector pairs, the extent of the overlapping region should be at least $4 \pi/M$ to ensure that it covers at least one injector pair in each sector. 
Choosing $\eta$ and the size of the overlapping region depends on the specific configuration, flow conditions, and number of detonation waves.
Based on our experiments in RDRE scenarios with a fixed number of co-rotating waves, an overlapping region covering two or four injectors is sufficient to ensure accurate approximations at interfaces. 
However, in scenarios with more complex dynamics, such as those involving subdominant or counter-propagating waves, a larger overlapping region may be necessary.
To determine the value of $\eta$ for a given extent of the overlapping region, one approach is to analyze the decay of the POD singular values as $\eta$ increases and
assess whether this leads to potential accuracy gains, keeping in mind that the computational cost of determining the regularization hyperparameters in the DD-$k$ OpInf learning problems grows exponentially with the number of sectors.

We first highlight the benefits of domain decomposition in terms of cost reduction and improved scalability.
Even though the training dataset comprises only $375$ snapshots, loading these snapshots into memory alone requires about $211$ GB of RAM.
Processing them and computing the POD basis via the thin SVD of the snapshot matrix requires even more RAM (we needed about $840$ GB in total).
This means that access to either a large shared-memory machine (for a serial or multi-threaded implementation) or a distributed-memory computer such as a cluster or a supercomputer (for a distributed implementation) is necessary to construct SD OpInf ROMs.
In contrast, by splitting the combustion chamber into four overlapping sectors, the full dataset can be processed in smaller and independent components in parallel.
The sector with the largest number of degrees of freedom requires roughly $64$ GB of RAM for loading the data into memory, which is a factor of $3.27 \times$ reduction compared to the single-domain approach.
This reduction decreases the overall memory and computing requirements and hence makes the construction of the domain-decomposed ROM more manageable.
We note that domain decomposition leads to memory and compute savings even if other methods such as randomized SVD~\cite{HMT11}, streaming SVD~\cite{Br06}, or the method of snapshots~\cite{BHL93} are used to compute the POD basis.

We next present the gains in terms of accuracy due to using domain decomposition.
The maximum reduced dimension for a quadratic SD OpInf ROM is $r=24$.
This value also represents the maximum reduced dimension for each sector in DD$-4$ OpInf with linear coupling terms~\eqref{eq:DD_OpInf_time_cont}.
These reduced dimensions respectively set the maximum number of operator coefficients that can be inferred via the OpInf regression problems~\eqref{eq:SD_OpInf_regularization_time_cont} and~\eqref{eq:DD_OpInf_subdom_regularization_time_cont}.
We therefore use $r \leq 24$ in SD OpInf and $r_1 = r_2 = r_3 = r_4 \leq 24$ in DD-$4$ OpInf in all our experiments.
\begin{remark}
As noted above, the coupling terms in our DD-$4$ OpInf experiments are linear; cf.~Eq.~\eqref{eq:DD_OpInf_time_cont}.
Considering quadratic coupling terms as well would double the size of the DD-$4$ OpInf ROM which in turn would further lower the maximum reduced dimensions to $r_1 = r_2 = r_3 = r_4 = 16$ due to the sparsity of the training dataset. 
Nevertheless, since the decomposition is done using overlapping, the respective POD bases contain information about the coupling between neighboring subdomains.
Therefore, as we will show in our results next, considering only linear coupling terms does not impact the accuracy of our DD-$4$ OpInf ROM.
\end{remark}

\begin{figure}[htb!]
\centering
\includegraphics[width=1.0\textwidth]{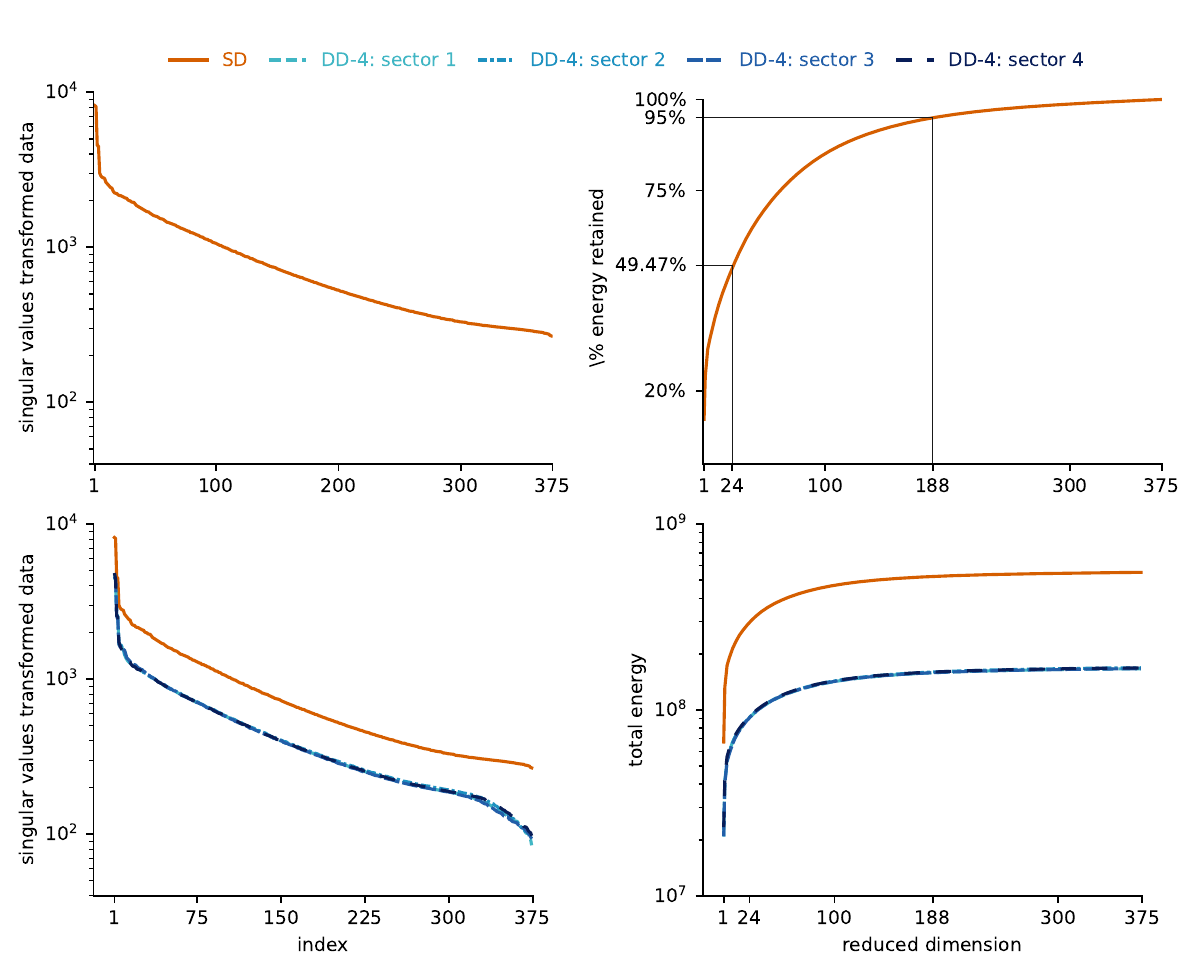}
\caption{POD singular values (top left) and retained energy (top right) for SD OpInf. POD singular values (bottom left) and total energy (bottom right) for both DD-$4$ and SD OpInf.}
\label{fig:RDE_SD_svals_and_ret_energy}
\end{figure}

Figure~\ref{fig:RDE_SD_svals_and_ret_energy} plots the POD singular values (top left) and corresponding retained energy (top right) for the single-domain approach.
It is not surprising that the POD singular values decay slowly, given the complexity of the underlying scenario and the large downsampling factor of the snapshots.
For example, a reduced dimension of $r = 188$ is necessary to ensure that $95\%$ of the total energy is retained.
In contrast, the maximum reduced dimension for SD OpInf ($r = 24$) retains only $49.47 \%$ of the total energy.
The POD singular values corresponding to the four sectors are plotted alongside the singular values for the single-domain approach in the bottom left figure.
The decomposition decreases their magnitude but does not lead to any significant decrease in the decay rate.
This is due to the symmetry of the domain and the (quasi-)periodicity of the corresponding solution fields. 
Lastly, the bottom right figure plots the corresponding total energy, which is consistent with what was observed in the singular value decay plot.
We note that even though the maximum reduced dimension retains only a small percentage of the total energy, the ROMs issue accurate predictions for quantities such as pressure and temperature, as we will show next.
This is because the POD basis is global, which means that it is constructed to represent all $n_s = 18$ transformed state variables~\eqref{eq:RDE_state_vars}.
Since the dynamics corresponding to variables such as pressure are smoother than the dynamics corresponding to species mass fractions, which are highly oscillatory, the global POD basis will represent some variables more accurately.
An alternative approach is to construct separate bases for individual variables or groups of variables, in which case the ROM dimension is given by the cumulative size of the respective bases.
This approach, however, is infeasible here due to the sparsity of the training set.

\begin{figure}[htb!]
\centering
\includegraphics[width=1.0\textwidth]{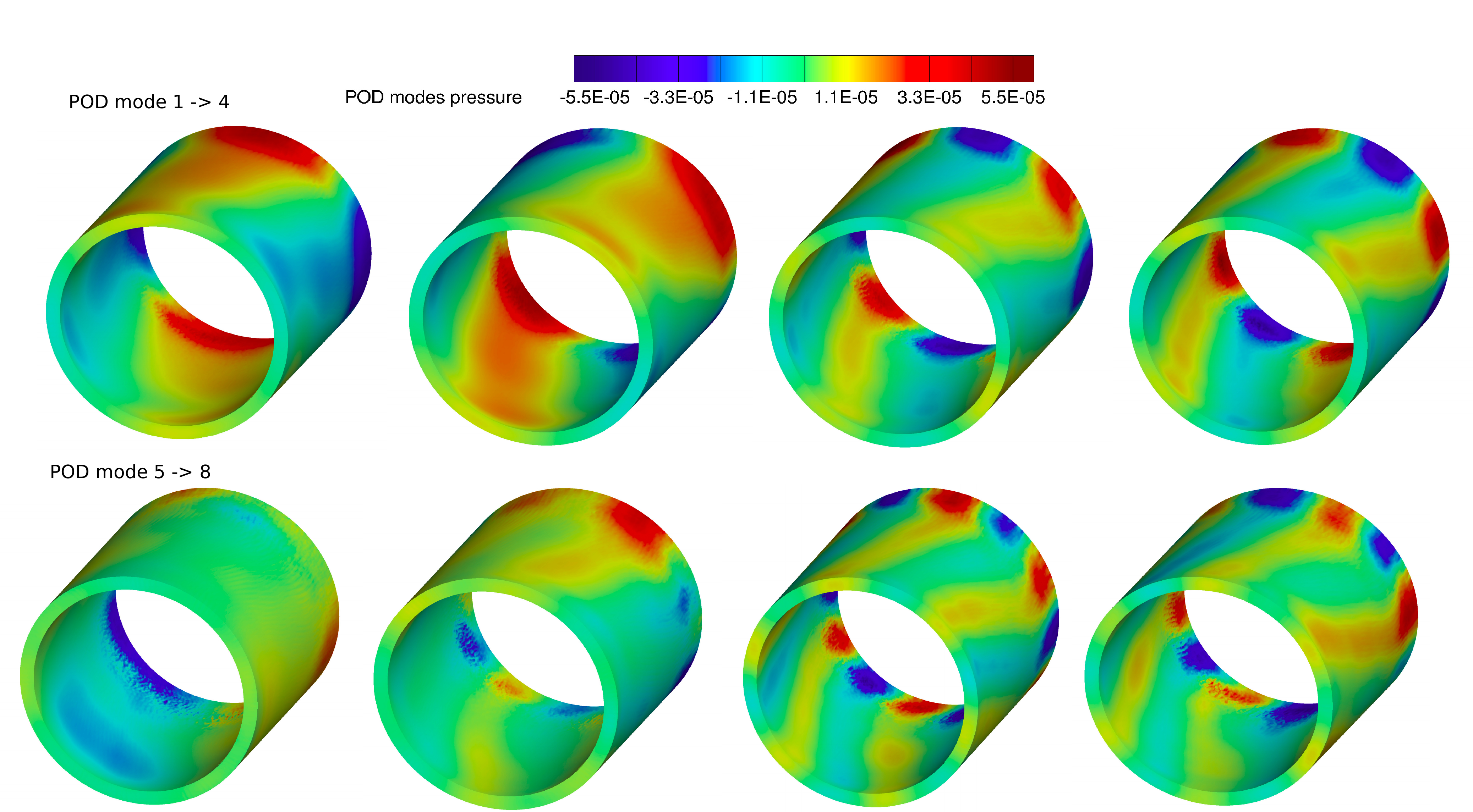}
\caption{First eight POD modes for the centered and scaled pressure field.}
\label{fig:RDE_first_eight_pressure_POD_modes}
\end{figure}

\begin{figure}[htb!]
\centering
\includegraphics[width=0.5\textwidth]{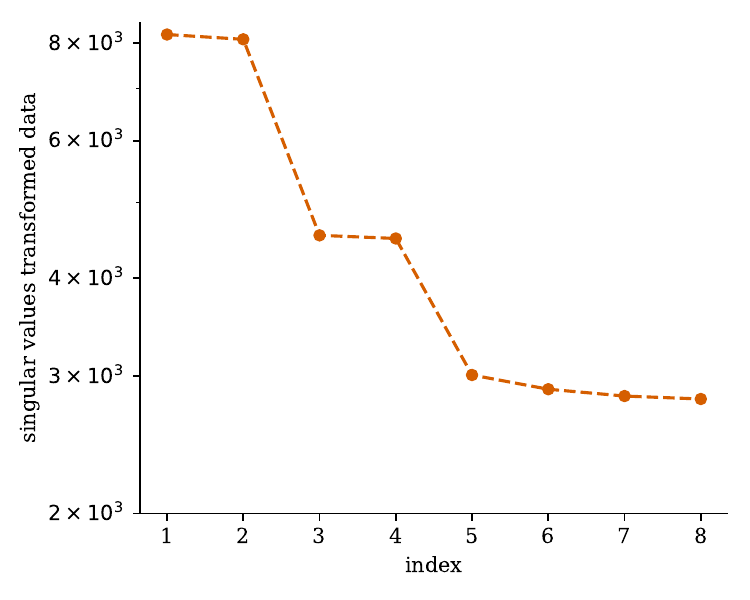}
\caption{First eight POD singular values for the centered and scaled pressure field.}
\label{fig:RDE_first_eight_POD_svals}
\end{figure}

\begin{figure}[htb!]
\centering
\includegraphics[width=1.0\textwidth]{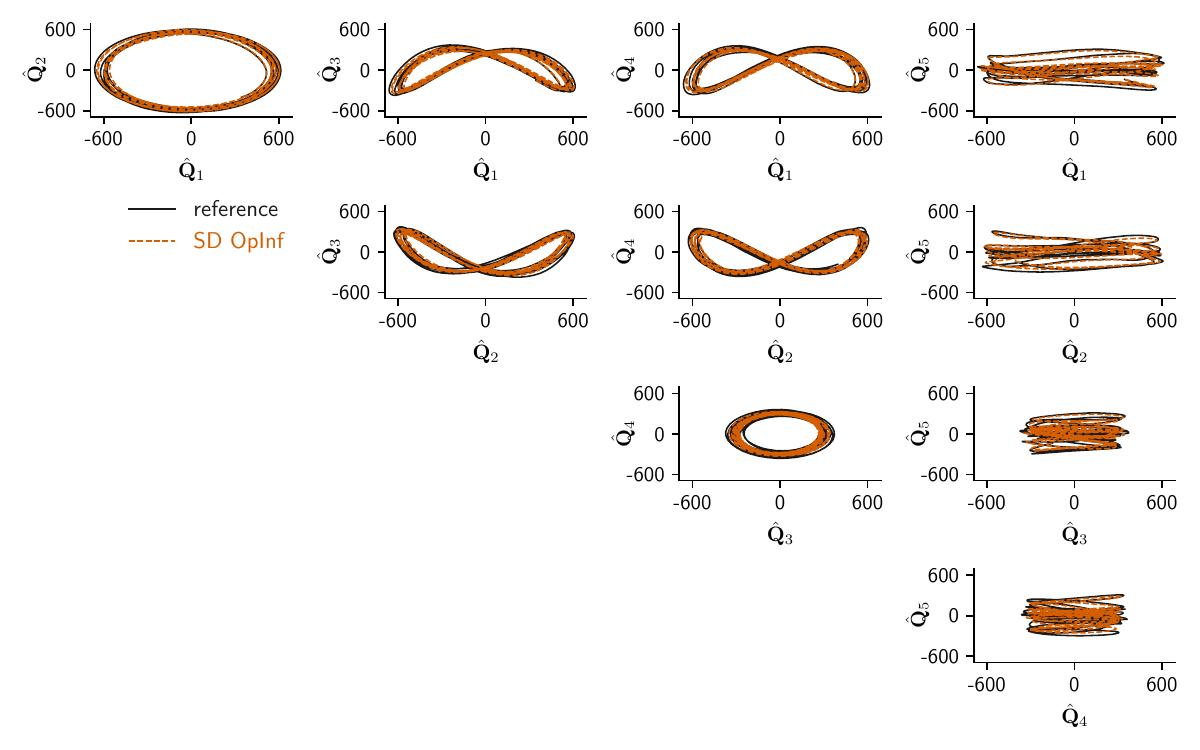}
\caption{Two-dimensional phase portraits of amplitudes over both training and prediction horizons corresponding to the five most energetic pressure POD modes.}
\label{fig:RDE_amplitudes_OpInf_2D}
\end{figure}

\begin{figure}[htb!]
\centering
\includegraphics[width=1.0\textwidth]{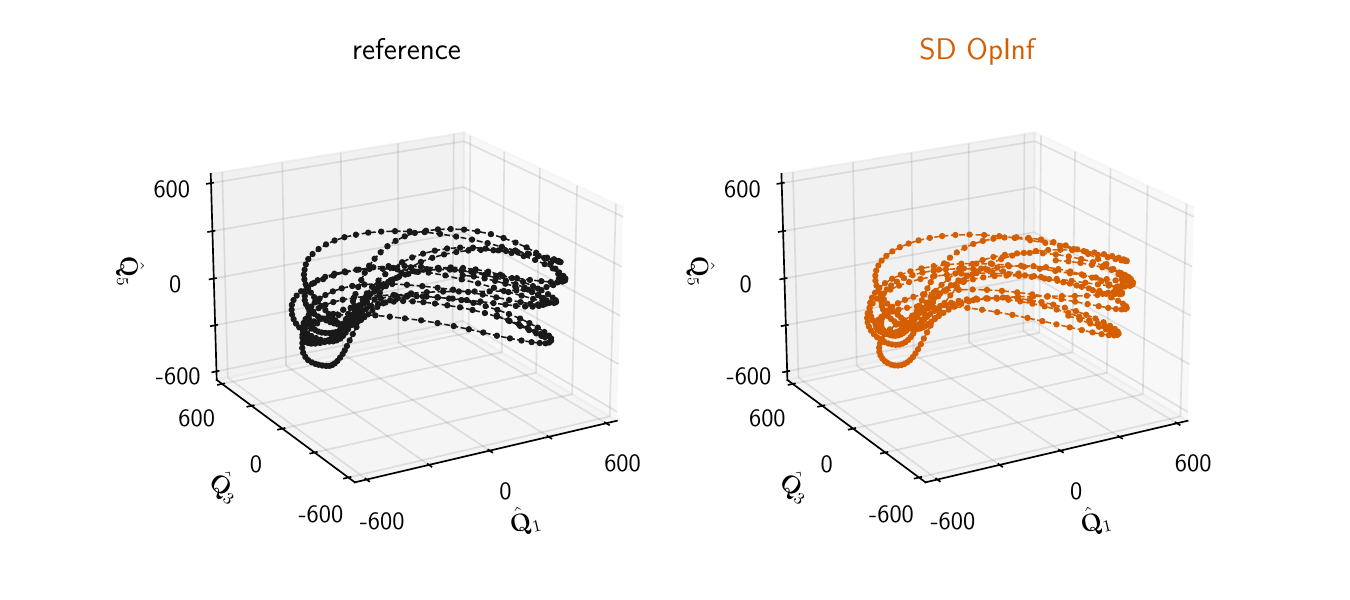}
\caption{Three-dimensional phase portrait of amplitudes for the first, third, and fifth most energetic pressure POD modes over both training and prediction horizons.}
\label{fig:RDE_amplitudes_OpInf_3D}
\end{figure}

Figure~\ref{fig:RDE_first_eight_pressure_POD_modes} plots the first eight POD modes for the centered (with respect to the mean over the training horizon) and scaled (with respect to the maximum absolute value of the centered data) pressure field; we refer the reader to Ref.~\cite{Ch21}, for example, for a discussion about pressure POD modes in the context of reduced modeling of swirl injectors.
We observe that modes $1$-$2$, $3$-$4$, and $7$-$8$ come in pairs, which is explained by the (quasi)-periodicity of the underlying dynamics.
In contrast, correlation between modes $5$ and $6$ seems poor, and this could be due to high amplitude fluctuations in coupled quantities such as entropy and enthalpy.
In general, correlation in phase space indicates that dynamical information captured in the respective modes is adequate to describe the dynamics of larger scale features. 
Lack of correlation indicates coupling and other forms of nonlinearity.
We furthermore plot, in Figure~\ref{fig:RDE_first_eight_POD_svals}, the first eight POD singular values corresponding to the SD OpInf approach.

We also compare the amplitudes --- over both training and prediction horizons --- obtained by projecting the centered and scaled snapshots via~\eqref{eq:project_snapshots} with those obtained via SD OpInf.
Figure~\ref{fig:RDE_amplitudes_OpInf_2D} plots two-dimensional phase portraits of the amplitudes corresponding to the five most energetic modes.
The pairs involving the first four amplitudes show strong correlations whereas the pairs involving the fifth amplitude are much less correlated.
This is consistent with the modes' structure observed in Figure~\ref{fig:RDE_first_eight_pressure_POD_modes}.
Furthermore, the SD OpInf approximate solutions closely match the reference results. 
While the 2D phase portraits involving the fifth amplitude do not exhibit significant correlations with other amplitudes, stronger correlations are observed in three-dimensional phase portraits. 
Figure~\ref{fig:RDE_amplitudes_OpInf_3D} plots the 3D phase portrait depending on the amplitudes corresponding to the first, third, and fifth most energetic modes.

\begin{figure}[htb!]
\centering
\includegraphics[width=1.0\textwidth]{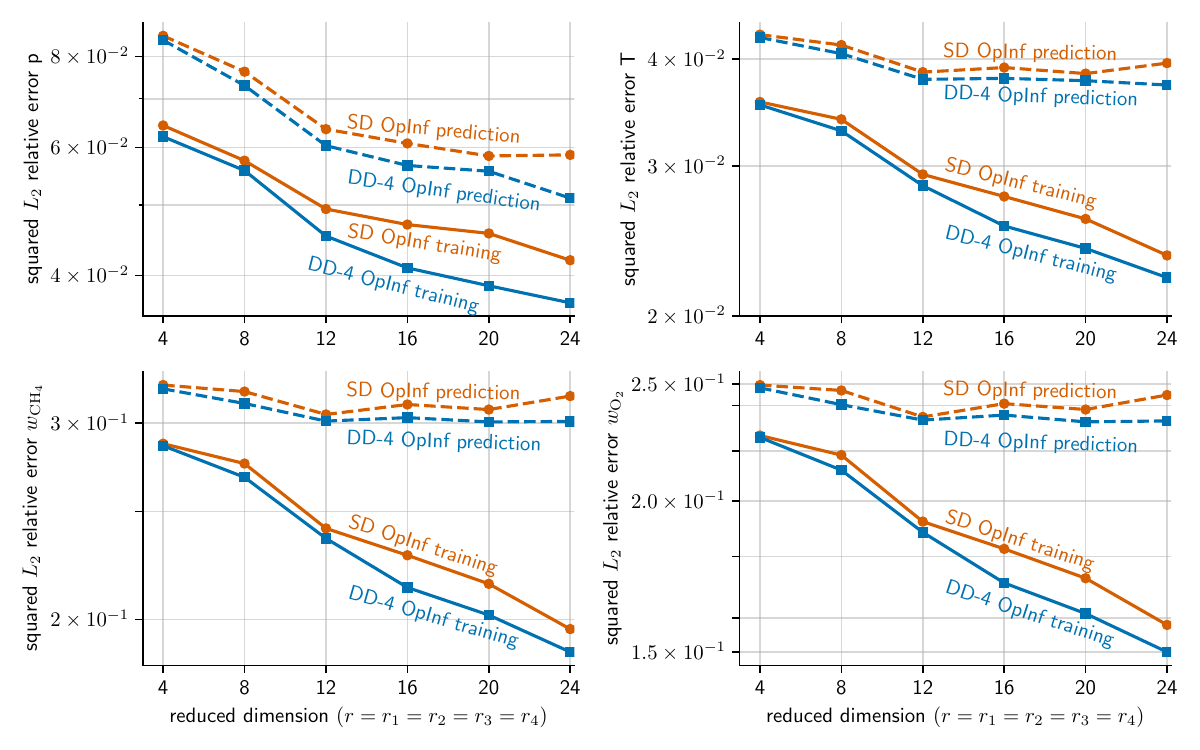}
\caption{Squared $L_2$ relative errors for pressure (top left), temperature (top right), and fuel (bottom left) and oxidizer (bottom right) mass fractions, for both training and prediction horizons.}
\label{fig:RDE_DD4_squared_rel_err}
\end{figure}

Figure~\ref{fig:RDE_DD4_squared_rel_err} plots the squared $L_2$ relative errors of the SD and DD-$4$ OpInf ROMs for pressure (top left), temperature (top right), and fuel (bottom left) and oxidizer (bottom right) mass fractions.
These four variables generally represent the most relevant quantities in RDRE simulations.
We note that these errors were computed with respect to the original state variables.
DD-$4$ OpInf leads to some accuracy gains for all four variables, for both training and prediction.
The squared relative error for pressure is decreased by up to $12.60 \%$ for training and up to $12.80 \%$ for prediction.
The gains for the other three variables are less significant, around $5\%$ for both training and prediction.
Nevertheless, these results indicate that domain decomposition can be used to construct accurate and scalable data-driven ROMs for complex and large-scale applications such as RDREs.

\begin{figure}[htb!]
\centering
\includegraphics[width=1.0\textwidth]{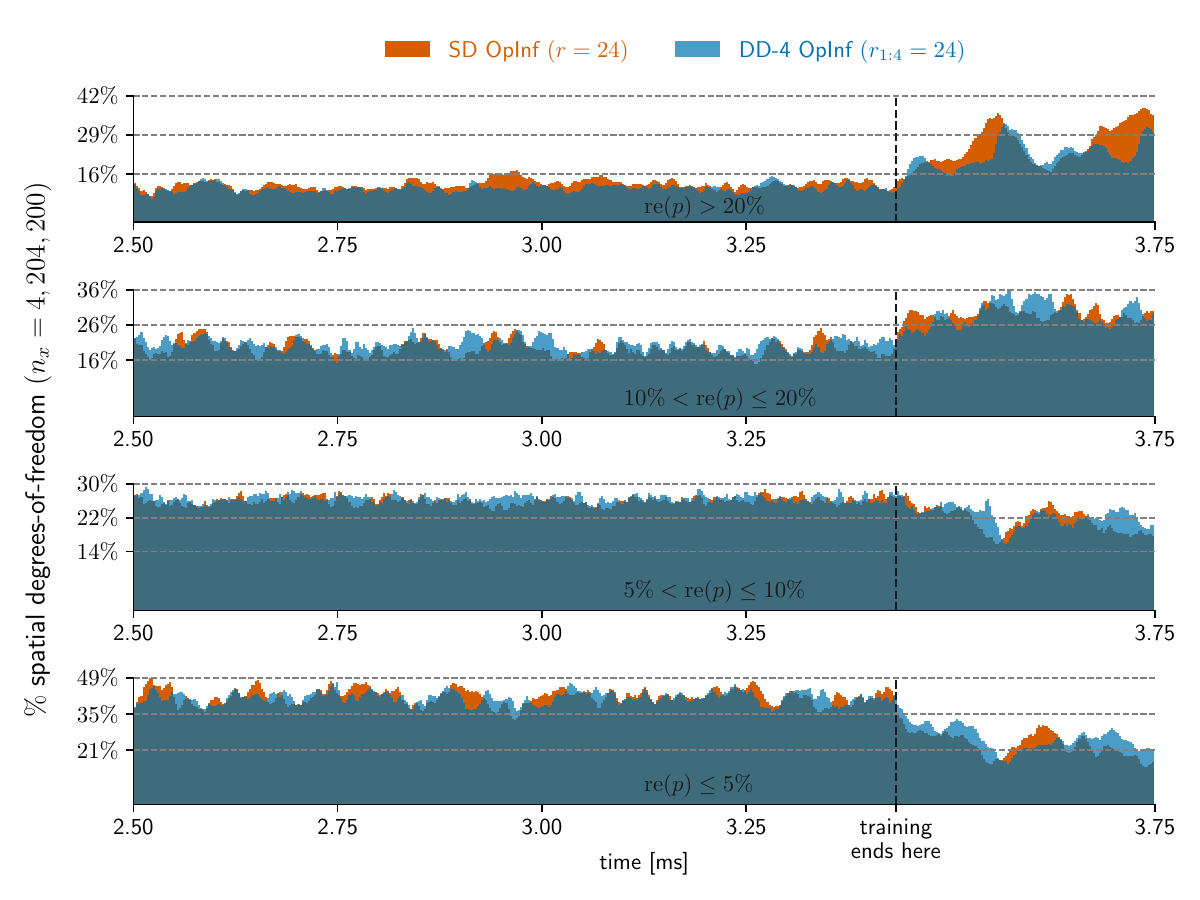}
\caption{Distribution of pointwise relative errors for pressure over time interval $[2.50, 3.75]$ ms. The dashed vertical lines mark the end of the training horizon.}
\label{fig:RDE_DD4_pressure_pointwise_rel_err}
\end{figure}

For a more detailed understanding of the accuracy of the SD and DD-$4$ OpInf ROMs, we also compute the distribution of pointwise relative errors in the entire spatial domain, over both training and prediction horizons.
Figure~\ref{fig:RDE_DD4_pressure_pointwise_rel_err} compares the spatio-temporal distribution of relative errors for pressure in the original coordinates, denoted by $\mathrm{re}(p)$.
From bottom to top, we plot the percentages of spatial degrees of freedom that satisfy $\mathrm{re}(p) \leq 5\%$ (bottom), $5\% < \mathrm{re}(p) \leq 10\%$ (second plot from the bottom), $10 < \mathrm{re}(p) \leq 20\%$ (second plot from the top), and $\mathrm{re}(p) > 20\%$ (top).
The percentages of spatial cells with large relative errors ($\mathrm{re}(p) > 20\%$) are decreased by up to $6.60\%$ (corresponding to $277,477$ spatial cells in total) over the training horizon and up to $15.20 \%$ over the prediction horizon (corresponding to $639,038$ spatial cells in total).
This, in turn, leads to an increase in the percentages of cells with small errors ($\mathrm{re}(p) \leq 10\%$).

\begin{figure}[htb!]
\centering
\includegraphics[width=1.0\textwidth]{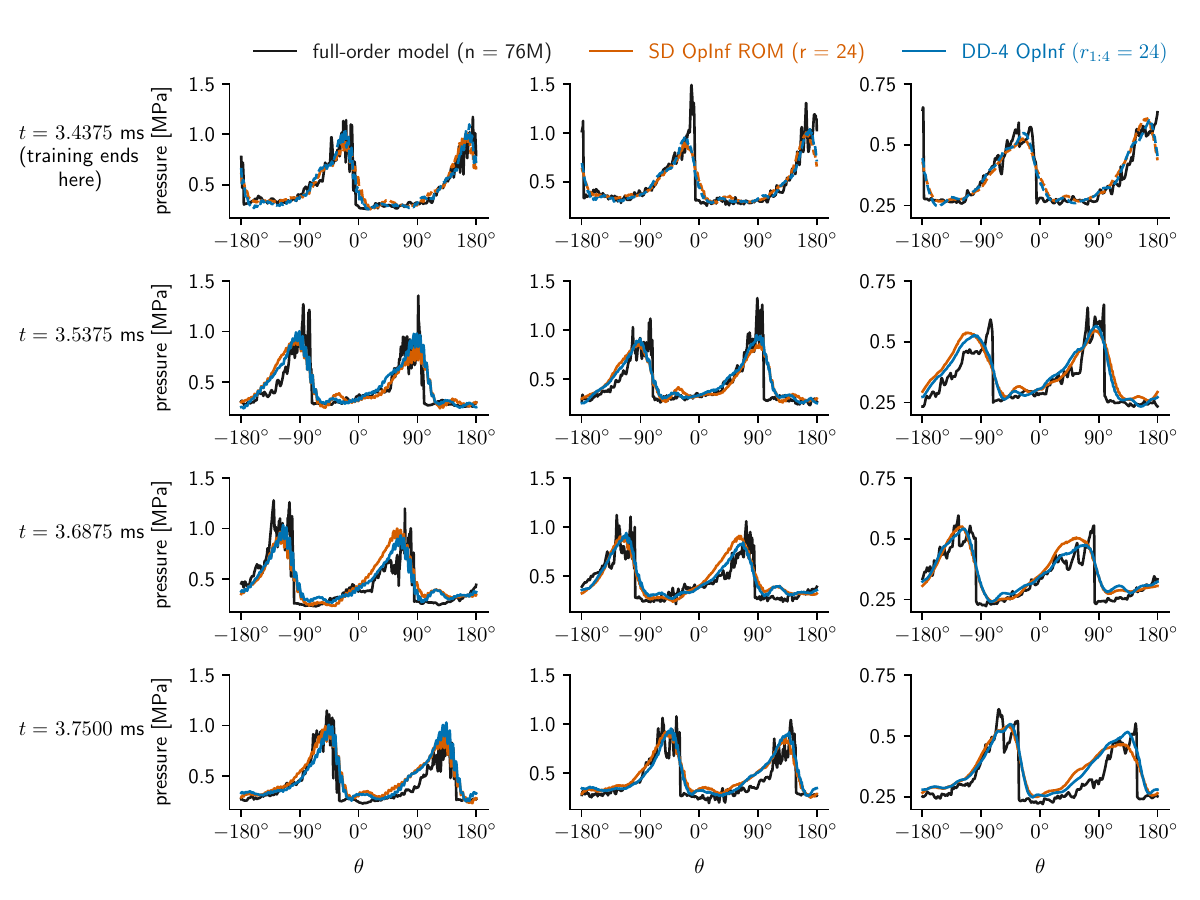}
\caption{One-dimensional circumferential pressure profile predictions. The columns plot the results at three representative locations close to the mid-channel. The rows plot the profiles at four representative time instants.}
\label{fig:RDE_DD4_pressure_profiles}
\end{figure}

\begin{figure}[htb!]
\centering
\includegraphics[width=1.0\textwidth]{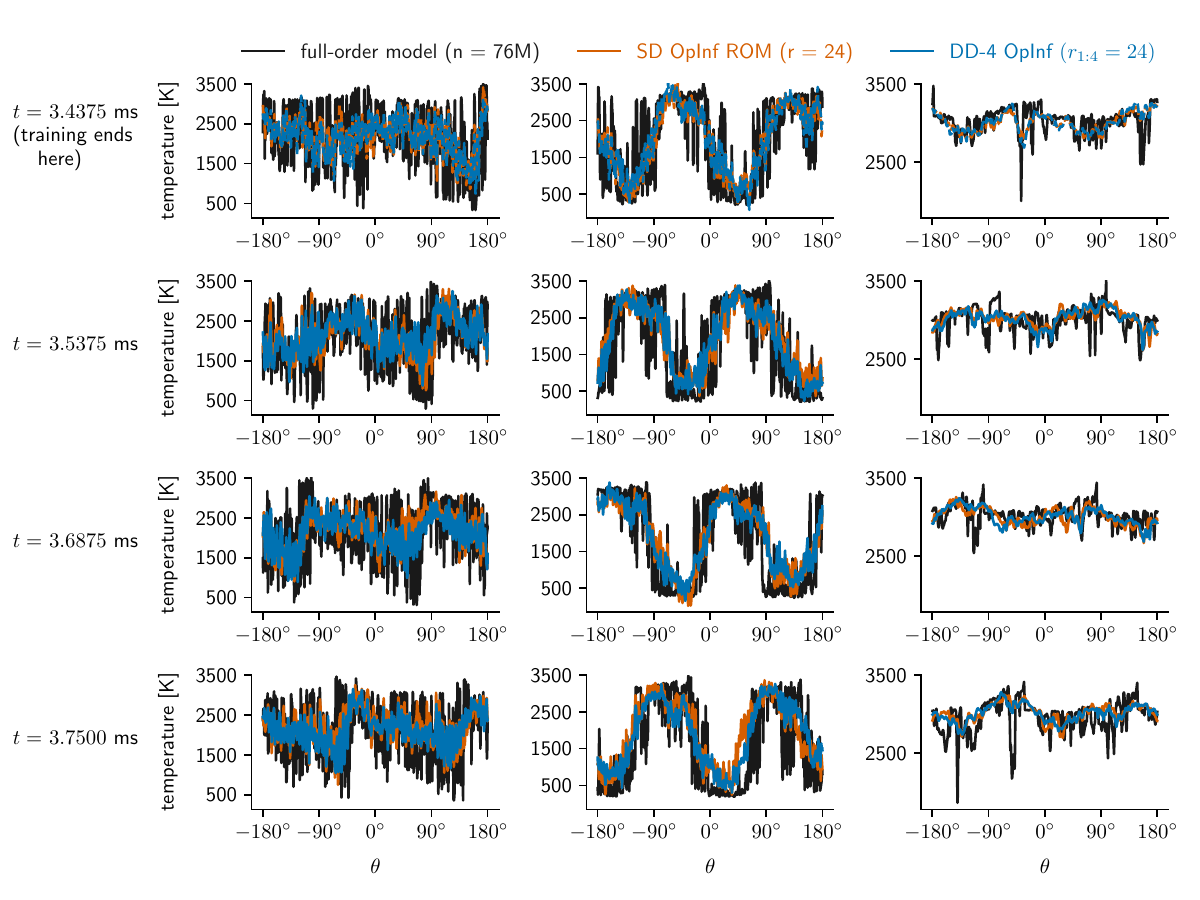}
\caption{One-dimensional circumferential temperature profile predictions. The columns plot the results at three representative locations close to the mid-channel. The rows plot the profiles at four representative time instants.}
\label{fig:RDE_DD4_temperature_profiles}
\end{figure}

\begin{figure}[htb!]
\centering
\includegraphics[width=1.0\textwidth]{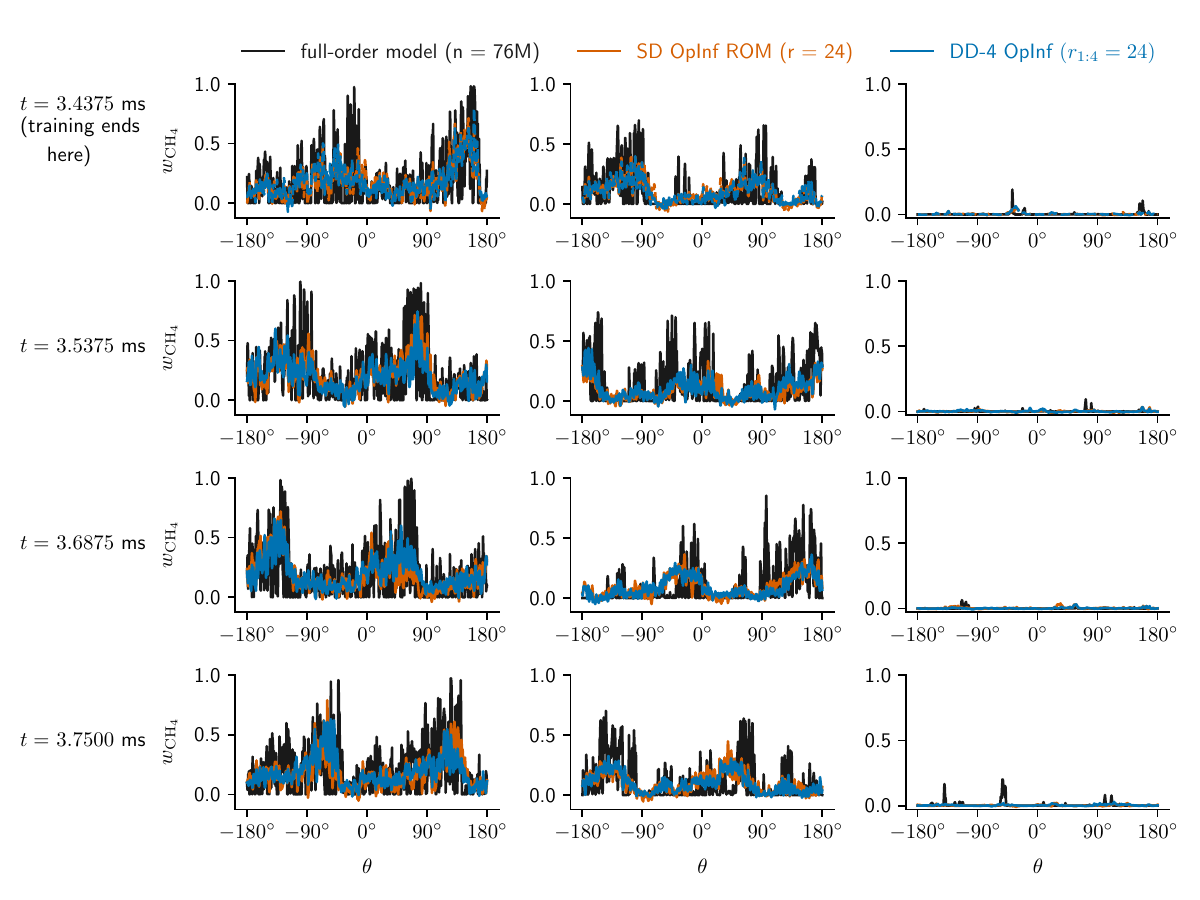}
\caption{One-dimensional circumferential $w_{\mathrm{CH_4}}$ profile predictions. The columns plot the results at three representative locations close to the mid-channel. The rows plot the profiles at four representative time instants.}
\label{fig:RDE_DD4_CH4_MF_profiles}
\end{figure}

\begin{figure}[htb!]
\centering
\includegraphics[width=1.0\textwidth]{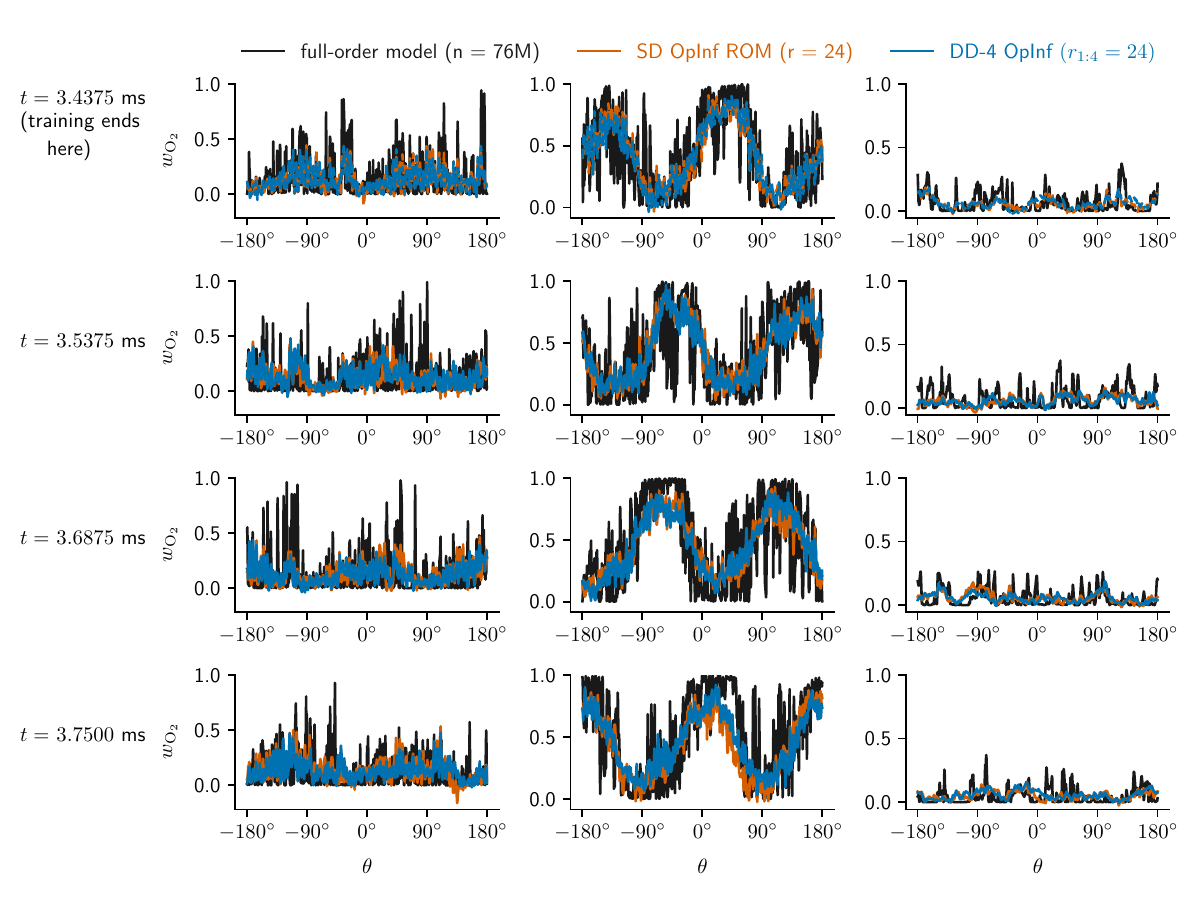}
\caption{One-dimensional circumferential $w_{\mathrm{O_2}}$ profile predictions. The columns plot the results at three representative locations close to the mid-channel. The rows plot the profiles at four representative time instants.}
\label{fig:RDE_DD4_O2_MF_profiles}
\end{figure}

We next use the ROM solutions obtained with the largest reduced dimensions ($r = r_1 = r_2 = r_3 = r_4 = 24$) to extract one-dimensional radial profiles at three representative locations close to the mid-channel. 
Figures~\ref{fig:RDE_DD4_pressure_profiles}--\ref{fig:RDE_DD4_O2_MF_profiles} plot these profiles for pressure (Figure~\ref{fig:RDE_DD4_pressure_profiles}), temperature (Figure~\ref{fig:RDE_DD4_temperature_profiles}), fuel mass fraction (Figure~\ref{fig:RDE_DD4_CH4_MF_profiles}), and oxidizer mass fraction (Figure~\ref{fig:RDE_DD4_O2_MF_profiles}).
Axially, the first location is close to the injectors, the second location is further way from the injectors but still within the detonation region, and the third location is downstream of the detonation zone.
The results corresponding to these three locations are respectively plotted in the first, second, and third column in each figure.
The four rows in each figure plot the profiles at four time instants.
The top row plots the profiles at the last time instant in the training horizon ($t = 3.5000$ ms).
The remaining three rows plot the profiles at three time instants in the prediction horizon: the $40$th ($t = 3.5375$ ms), the $100$th ($t = 3.6875$ ms), and the $126$th ($t = 3.7500$ ms) time instant.
For pressure, the DD-$4$ OpInf accuracy gains observed in Figures~\ref{fig:RDE_DD4_squared_rel_err} and~\ref{fig:RDE_DD4_pressure_pointwise_rel_err} correspond to more accurate predictions for both the frequency and amplitude of the profiles. 
For the other three quantities, we do not observe a significant gain in accuracy due to DD-$4$ OpInf, which is consistent with the error plots in Figures~\ref{fig:RDE_DD4_squared_rel_err}. 
Overall, the ROMs perform well given the complexity of the underlying dynamics, which is clearly visible in the plotted profiles, especially at the location closer to the injectors.

Lastly, the runtimes of the SD and DD-$4$ OpInf ROMs, averaged over $1,000$ runs, are $0.0296 \pm 5.012 \times 10^{-4}$ and $0.1256 \pm 4.4465 \times 10^{-4}$ seconds respectively.
This translates into $12$ (for SD OpInf) and $11$ (for DD-$4$ OpInf) orders of magnitude reduction compared to the runtime of the high-fidelity model.

\section{Conclusions} \label{sec:conclusions}
This paper introduced a domain decomposition formulation into the construction of a data-driven reduced model of large-scale systems.
Our results for a large-scale simulation of a three-dimensional unsteady rotating detonation rocket engine combustion chamber with a sparse training dataset showed that domain decomposition improves the accuracy of the equivalent single-domain approach while decreasing the overall computational burden for preprocessing the training data. 
This reduction in computational costs is essential for enabling the construction of data-driven reduced models of large-scale applications.
Overall, given the complexity of the scenario under consideration, our results show that domain decomposition can indeed provide a viable approach to constructing accurate, predictive, and scalable data-driven reduced models of large-scale systems.
While these results illustrate the potential of the approach, they also highlight a number of open questions.
Important future directions of work include alternative formulations of the domain-decomposed problem that explicitly enforce continuity across subdomains using constraints, strategies to optimally decompose the domain in a data-driven adaptive fashion, and an approach to identify optimal data scalings and transformations that may embed physical constraints in the reduced-order basis.   
In the broader context of aerospace propulsion, constructing reduced models that maintain the desired engineering accuracy while drastically decreasing the cost of computationally expensive high-fidelity simulations can have a strong impact and ultimately enable tasks such as design optimization which would be almost impossible to perform otherwise.

\section*{Acknowledgments}
This work was supported in part by AFRL Grant FA9300-22-1-0001 and the Air Force Center of Excellence on Multifidelity Modeling of Rocket Combustor Dynamics under grant FA9550-17-1-0195.
The views expressed are those of the author and do not necessarily reflect the official policy or position of the Department of the Air Force, the Department of Defense, or the U.S. government.

Distribution Statement A: Approved for Public Release; Distribution is Unlimited. PA$\#$ AFRL-2023-5401

\clearpage
\bibliography{DD_OpInf_final}

\end{document}

%% file: settings.tex
\usepackage[utf8]{inputenc}
\usepackage{textcomp}
\usepackage{tikz}
\usetikzlibrary{shapes,arrows}
\usepackage{graphicx}
\usepackage{amsmath}
\usepackage[version=4]{mhchem}
\usepackage{siunitx}
\usepackage{combelow}
\usepackage{multirow}
\usepackage{multicol}
\usepackage{xcolor}
\usepackage{tikzpagenodes}
\usepackage{amsfonts}
\usepackage{todonotes}
\usepackage{mathtools}
\usepackage{arydshln}
\usepackage{algorithm}
\usepackage{algpseudocode}
\usepackage{nomencl}
\usepackage{hyperref}
\usepackage{fancyhdr}
\usepackage{longtable,tabularx}
\newtheorem{remark}{Remark}

\newcommand{\argmin}{\mathop{\mathrm{argmin}}\limits}

\hypersetup{
    pdftitle={Domain decomposition for data-driven reduced modeling of large-scale systems},
    pdfauthor={Ionut-Gabriel Farcas, Rayomand Gundevia, Ramakanth Munipalli, Karen Willcox},
    pdfkeywords={domain decomposition, domain-decomposed data-driven reduced modeling, data-driven reduced modeling, scientific machine learning, large-scale simulation, rotating detonation rocket engines, combustion}
}
